\documentclass[12pt, reqno]{amsart}
\usepackage[all]{xy}
\usepackage{amsmath}
\usepackage{amssymb}
\usepackage{amsfonts}
\usepackage{amsthm}
\usepackage{amscd}
\usepackage{bm}
\usepackage{ulem}
\usepackage{graphicx}
\usepackage[dvipdfmx]{color}
\newtheorem{thrm}{Theorem}

\newtheorem{prop}{Proposition}

\newtheorem{thm}{Theorem}[section]
\newtheorem{lem}[thm]{Lemma}
\newtheorem{pro}[thm]{Proposition}

\newtheorem{rem}[thm]{Remark}

\newtheorem{defn}[thm]{Definition}
\newtheorem{prob}[thm]{Problem}

\newcommand{\Z}{\mathbf{Z}}
\newcommand{\Q}{\mathbf{Q}}
\newcommand{\N}{\mathbf{N}}

\begin{document}
\title[the special derivation Lie algebras of the free Lie algebra]
{On the abelianization of the special derivation \\ Lie algebras of free Lie algebras}
\author{Naoya Enomoto}
\address{Naoya Enomoto; The University of Electro-Communications,
1-5-1, Chofugaoka, Chofu city, Tokyo 182-8585, Japan.}
\email{enomoto-naoya@uec.ac.jp}
\author{Takao Satoh}
\address{Takao Satoh; Department of Mathematics, Faculty of Science Division II, Tokyo University of Science,
1-3, Kagurazaka, Shinjuku-ku, Tokyo 162-8601, Japan.}
\email{takao@rs.tus.ac.jp}
\subjclass[2020]{17B40(Primary), 20F14, 20C30, 17B56(Secondary)}%ams classifcation
\keywords{Derivation algebras, Free Lie algebras.}
\maketitle
\begin{abstract}
In this paper, we show that there are infinitely many linearly independent elements
in the abelianization of the Lie algebra of special derivations of a free Lie algebra by using the Morita traces.
Furthermore, we show that the abelianization contains non-trivial elements which are killed by the Morita traces. 
\end{abstract}
\tableofcontents
\section{Introduction}

Over the last forty years, the Johnson homomorphisms have been studied for the mapping class groups of surfaces and automorphism groups of free groups
by many authors over generations.
To put it simply, the Johnson homomorphisms describe a series of approximations of the graded quotients of a certain
filtrations of the groups. Since each of the Johnson homomorphisms is injective by definition,
to determine their images and cokernels is a natural problem.
In order to attack the problem, Morita \cite{Mo1} originally introduced trace maps.
Morita's traces played a crucial role in the study of the Johnson homomorphisms, and gave many contributions to the development of the field.
(For surveys of the field, see Morita \cite{Mo3}, Hain \cite{Ha2} and Satoh \cite{Sa1} for example.)

\vspace{0.5em}

Let us fix notation.
Let $F_n$ be the free group of rank $n$ generated by $x_1, \ldots, x_n$, and $\mathrm{Aut}\,F_n$  the automorphism group of $F_n$.
Set $H:=H_1(F_n,\Z)$, the abelianization of $F_n$. Denote by $\mathcal{L}_n=\bigoplus_{k \geq 1} \mathcal{L}_n(k)$ the free Lie algebra generated by $H$.
We have a certain filtration
$\mathcal{A}_n(1) \supset \mathcal{A}_n(2) \supset \cdots$ of $\mathrm{Aut}\,F_n$, so-called the Andreadakis-Johnson filtration of $\mathrm{Aut}\,F_n$,
and for each $k \geq 1$ we have the $k$-th Johnson homomorphism
\[ \tau_k : \mathrm{gr}^k(\mathcal{A}_n):= \mathcal{A}_n(k)/\mathcal{A}_n(k+1) \rightarrow \mathrm{Hom}_{\Z}(H, \mathcal{L}_n(k+1)). \]
Furthermore, the graded sum $\bigoplus_{k \geq 1} \mathrm{gr}^k(\mathcal{A}_n) \rightarrow \mathrm{Der}^+(\mathcal{L}_n)$ can be considered as
a Lie algebra homomorphism, where $\mathrm{Der}^+(\mathcal{L}_n)$ is the positive part of the derivation Lie algebra of $\mathcal{L}_n$.
Let 
\[ \mathcal{C}_n(k) := H^{\otimes k}/ \langle a_1 \otimes a_2 \cdots \otimes a_k - a_k \otimes a_1 \cdots \otimes a_{k-1} \,|\, a_i \in H \rangle \]
be the quotient of $H^{\otimes k}$ by the permutation action of the cyclic group $C_k$ of degree $k$ on the components.
By combining results by Satoh \cite{S03, S13} and Darn\'{e} \cite{Dar}, we see that the cokernel of $\tau_k$ is isomorphic to
$\mathcal{C}_n(k)$ for $n \geq k+2$, which is detected by the trace map
\[ \mathrm{Tr}_k : H^* \otimes_{\Z} \mathcal{L}_n(k+1) \rightarrow H^* \otimes_{\Z} H^{\otimes (k+1)} \rightarrow H^{\otimes k} \rightarrow \mathcal{C}_n(k) \]
defined as the compositions of maps including a contraction.
(See Section {\rmfamily \ref{S-John-trace}} for details.)

\vspace{0.5em}

Let $S^k H$ be the symmetric power of $H$ of degree $k \geq 1$.
The composition of $\mathrm{Tr}_k$ and the natural projection $\mathcal{C}_n(k) \rightarrow S^kH$, is called the Morita trace as mentioned above,
denoted by $\mathrm{MT}_k : H^* \otimes_{\Z} \mathcal{L}_n(k+1) \rightarrow S^k H$.
Morita tried to give one of characterizations of the Morita traces as generators of the abelianization of $\mathrm{Der}^+(\mathcal{L}_n)$ as a Lie algebra.
In \cite{Mo3}, Morita constructed a surjective Lie algebra homomorphism
\[ \Theta := \mathrm{id}_1 \oplus \bigoplus_{k \geq 2} \mathrm{MT}_k : H_1(\mathrm{Der}^+(\mathcal{L}_n),\Z) \rightarrow \mathrm{Im}(\tau_1) \oplus
   \bigoplus_{k \geq 2} S^k H \]
where $\mathrm{id}_1$ is the identity map on the degree one part, and the target is understood to be an abelian Lie algebra.
Morita \cite{Mo3} conjectured that the above map induces an isomorphism,
and showed that it is true up to degree $n(n-1)$ over $\Q$ by using a work of Kassabov in \cite{Kas}.
Morita-Sakasai-Suzuki \cite{MSS} showed that $\Theta$ induces an isomorphism up to degree $n-2$ over $\Z$.
We remark that for the free metabelian Lie algebra $\mathcal{L}_n^M$, Enomoto-Satoh \cite{ES1} showed that for any $n \geq 4$
\[ H_1(\mathrm{Der}^+(\mathcal{L}_n^M),\Z) \cong \mathrm{Im}(\tau_1) \oplus \bigoplus_{k \geq 2} S^k H. \]

For a subgroup $\mathcal{S} \leq \mathrm{Aut}\,F_n$, we can consider subgroups $\mathcal{S} \cap \mathcal{A}_n(k)$, and the $k$-th
Johnson homomorphism
\[ \mathrm{gr}^k(\mathcal{S}) := \mathcal{S} \cap \mathcal{A}_n(k)/\mathcal{S} \cap \mathcal{A}_n(k+1)
  \rightarrow \mathrm{Hom}_{\Z}(H, \mathcal{L}_n(k+1)) \]
by the same argument as that of $\tau_k$.
Let $\mathrm{LB}_n$ be the mapping class group of a collection $C$ of $n$ disjoint,
unknotted, oriented circles embedded in the three ball $B^3$. 
Namely, $\mathrm{LB}_n$ is the group of isotopy classes of orientation preserving homeomorphisms of $B^3$ which fix $C$ setwise.
The group $\mathrm{LB}_n$ is called the loop braid group on $n$ components.
By a work of Goldsmith \cite{Gol}, $\mathrm{LB}_n$ can be embedded into $\mathrm{Aut}\,F_n$, and its image is given by
\[ \{ \sigma \in \mathrm{Aut}\,F_n \,|\, x_i^{\sigma} = a_i x_{\mu(i)} a_i^{-1},
      \hspace{0.5em} a_i \in F_n, \hspace{0.5em} \mu \in \mathfrak{S}_n \hspace{0.5em} (1 \leq i \leq n) \}. \]
Through this embedding, we consider $\mathrm{LB}_n$ as a subgroup of $\mathrm{Aut}\,F_n$.
We remark that Fenn-Rim\'{a}nyi-Rourke \cite{FRR} obtained a finite presentation of $\mathrm{LB}_n$, and
called $\mathrm{LB}_n$ the braid-permutation group. (See \cite{Dam} for details.)
As is well-known, Artin's braid group $B_n$ of $n$ strands can be embedded into $\mathrm{Aut}\,F_n$.
If we consider $B_n$ as a subgroup of $\mathrm{Aut}\,F_n$, then we have $B_n \subset \mathrm{LB}_n \subset \mathrm{Aut}\,F_n$.

\vspace{0.5em}

In our previous paper \cite{ES3}, we studied the Johnson homomorphisms of $\mathrm{LB}_n$, denoted by $\tau_k^P$.
The graded sum of the images of $\tau_k^P$ is contained in a Lie subalgebra of $\mathrm{Der}^+(\mathcal{L}_n)$,
called the tangential derivation Lie algebra of $\mathcal{L}_n$, denoted by $\mathfrak{p}_n=\bigoplus_{k \geq 1} \mathfrak{p}_n(k)$.
Alekseev and Torossian \cite{AT1, AT2} studied $\mathfrak{p}_n$ in a series of their works including for the Kashiwara-Vergne conjecture, and
showed $\mathfrak{p}_n$ is a Lie subalgebra of $\mathrm{Der}^+(\mathcal{L}_n)$.
In a different context, Satoh \cite{S18} considered $\mathfrak{p}_n$ for the study of the Johnson homomorphisms of McCool groups.
In \cite{ES3}, we introduced a certain trace maps $\widetilde{\mathrm{Tr}}_k$ such that
$\mathrm{Im}(\tau_k^P) \subset \mathrm{Ker}(\widetilde{\mathrm{Tr}}_k) \subset \mathfrak{p}_n(k)$ for $k \geq 2$.
By using these trace maps, we determined the structure of the images and the cokernels of
$\tau_k^P$ for $1 \leq k \leq 4$. At the present stage, however, the Johnson cokernels of $\mathrm{LB}_n$ have not determined yet in general.
%We will give some descriptions for the abelianization of $\mathfrak{p}_n$ in a forthcoming paper \cite{ES5}.

\vspace{0.5em}

As a subsequent study, this paper considers the Johnson homomorphisms of Artin's braid groups $B_n$, denoted by $\tau_k^B$.
The graded sum of the images of $\tau_k^B$ is contained in a Lie subalgebra of $\mathfrak{p}_n$,
called the special derivation Lie algebra of $\mathcal{L}_n$, denoted by $\mathfrak{b}_n=\bigoplus_{k \geq 1} \mathfrak{b}_n(k)$.
The Lie algebra $\mathfrak{b}_n$ seemed to be introduced by Ihara in a series of his researches for the absolute Galois group over $\Q$ in the 1980s
In particular, Ihara \cite{Iha} calls elements of $\mathfrak{b}_n$ normalized special derivations.
We determined the images and the cokernels of
$\tau_k^B: \mathrm{gr}^k(B_n) \rightarrow \mathfrak{b}_n(k)$. 
By combining a work of Darn\'{e} \cite{Da2} and independent works of Kohno \cite{Koh} and Falk-Randell \cite{F-R},
we can calculate the rank of a free abelian group $\mathrm{gr}^k(B_n)$ for each $k \geq 1$. (See Section {\rmfamily \ref{S-M-PB} for details.)
Darn\'{e} \cite{Da3} showed that $\tau_1^B$ and $\tau_2^B$ are surjective, and gave a description of the ranks of the Johnson cokernels. 
In this paper, we describe the cokernel of $\tau_3^B$ and $\tau_4^B$ by using trace maps.
\begin{prop}[$=$ Propositions {\rmfamily \ref{T-John-1}} and {\rmfamily \ref{T-John-2}}]\label{T-I-1}
The sequences
\[ 0 \rightarrow \mathrm{gr}^3(B_n) \xrightarrow{\tau_3^B} \mathfrak{b}_n(3) \xrightarrow{\mathrm{MT}_3^B} \mathrm{Im}(\mathrm{MT}_3^B) \rightarrow 0 \]
for any $n \geq 3$ and
\[ 0 \rightarrow \mathrm{gr}_{\Q}^4(B_n) \xrightarrow{\tau_{4,\Q}^B} \mathfrak{b}_n^{\Q}(4) \xrightarrow{\mathrm{Tr}_{[2,1^2]}^{\Q}}
    \mathrm{Im}(\mathrm{Tr}_{[2,1^2]}^{\Q}|_{\mathfrak{b}_n^{\Q}(4)}) \rightarrow 0 \]
for any $n \geq 4$ are $\mathfrak{S}_n$-equivariant exact sequences.
\end{prop}
\noindent
Here the subscript and the superscript $\Q$ means tensoring with $\Q$ over $\Z$, and $\mathrm{Tr}_{[2,1^2]}$ is a kind of trace map, for the definition, see (\ref{eq-def-tr_4}) below.

In addition to show the above proposition, in order to investigate the abelianization of $\mathfrak{b}_n$,
we need a basis of each of $\mathfrak{b}_n(k)$.
In Section A.1 in \cite{Da3}, Darn\'{e} gave a formula for the rank of $\mathfrak{b}_n(k)$ for each $k \geq 1$, which
recovers a classical formula for Milnor invariants by Orr. (See Theorem 15 in \cite{Orr}.)
In general, in contrast to the case of the degree $k$ part of $\mathrm{Der}^+(\mathcal{L}_n)$, to write down generators
of $\mathfrak{b}_n(k)$ explicitly is quite complicated even in the low degree $k$.
We consider to decompose $\mathfrak{b}_n(k)$ into $\mathfrak{S}_n$-submodules with compositions of $k+1$.
In Section {\rmfamily \ref{S-strB}}, we give a basis of each component $\mathfrak{b}_n(k, \alpha)$ of $\mathfrak{b}_n(k)$
for any $1 \leq k \leq 4$ and any composition $\alpha$.
As for the Morita trace map $\mathrm{MT}_{k}^B$, we show the following.
\begin{thrm}[$=$ Theorem {\rmfamily \ref{T-im-mtr}}]\label{T-I-2}
For any $n \geq 3$ and any odd $k \geq 3$, we have
\[ \mathrm{rank}_{\Z}(\mathrm{Im}(\mathrm{MT}_{k}^B)) = \binom{n+k-1}{k} - \dfrac{n(n+1)}{2}. \]
\end{thrm}

\vspace{0.5em}

In the latter part of the paper, we consider the abelianization of the Lie algebra $\mathfrak{b}_n$.
For each $k \geq 1$, let $H_1(\mathfrak{b}_n,\Z)_k$ be the degree $k$ part of $H_1(\mathfrak{b}_n,\Z)$.
An immediate consequence obtained from Theorem {\rmfamily \ref{T-I-2}}, by considering the restriction of $\Theta$ to $\mathfrak{b}_n$,
we see that $H_1(\mathfrak{b}_n,\Z)_k$ contains non-trivial free abelian group of finite rank, and that
$H_1(\mathfrak{b}_n,\Z)$ is not finitely generated as a $\Z$-module. (For details, see Theorem {\rmfamily \ref{T-strB-2}}.)
We should remark that Yusuke Kuno communicated to us the non-finite generation of $H_1(\mathfrak{b}_2,\Z)$ can be obtained with Soul\'{e} elements defined by Ihara \cite{Iha}.
(See Remark {\rmfamily \ref{R-kuno}} for details.)

\vspace{0.5em}

In Subsection {\rmfamily \ref{Ss-ab-low}}, we show that for any $2 \leq k \leq 4$, $H_1(\mathfrak{b}_n,\Z)_k$ is trivial
except for the component that can be detected by the Morita trace $\mathrm{MT}_3^B$.
In Subsection {\rmfamily \ref{Ss-ab-5}}, we show that there exist linearly independent elements in $\mathrm{Ker}(\mathrm{MT}_5^B)$ and
not in $[\mathfrak{b}_n(4),\mathfrak{b}_n(1)]$.
\begin{thrm}\label{T-I-3}($=$ Theorem {\rmfamily \ref{T-ninj-MT}})
For $n \geq 3$, the surjective Lie algebra homomorphism
\[ \overline{\Theta|_{\mathfrak{b}_n}} : H_1(\mathfrak{b}_n, \Z) \rightarrow \mathrm{Im}(\tau_1) \oplus
   \bigoplus_{m \geq 1} \mathrm{Im}(\mathrm{MT}_{2m+1}^B)  \]
induced from $\Theta|_{\mathfrak{b}_n}$ is not injective.
\end{thrm}
\noindent
This implies that the structure of the abelianization of $\mathfrak{b}_n$ is much more complicated than that of $\mathrm{Der}^+(\mathcal{L}_n^M)$
and $\mathrm{Der}^+(\mathcal{L}_n)$. We expect that $\mathrm{Ker}(\mathrm{MT}_{5,\Q}^B)/[\mathfrak{b}_n^{\Q}(4), \mathfrak{b}_n^{\Q}(1)]$ is of dimension
$\binom{n}{3}+\binom{n}{4}$. (For details, see Remark {\rmfamily \ref{R-ab-5-5}}.)

\tableofcontents

\section{Notation and conventions}\label{S-Not}

In this section, we fix notation and conventions.
Throughout this paper, a group means a multiplicative group if otherwise noted.
Let $G$ be a group. 

\begin{itemize}
\item The abelianization of $G$ is denoted by $G^{\mathrm{ab}}$ if otherwise noted.
\item The automorphism group $\mathrm{Aut}\,G$ acts on $G$ from the right. For any $\sigma \in \mathrm{Aut}\,G$ and $x \in G$,
      the action of $\sigma$ on $x$ is denoted by $x^{\sigma}$.
\item For a normal subgroup $N$, we often denote the coset class of an element $g \in G$ by the same $g$ in the quotient group $G/N$ if
      there is no confusion.
\item For elements $x$ and $y$ of $G$, the commutator bracket $[x,y]$ of $x$ and $y$
      is defined to be $[x,y]:=xyx^{-1}y^{-1}$. Then for any $x, y, z \in G$, we have
\begin{equation}
  [xy,z] =[x,[y,z]][y,z][x,z], \hspace{1em} [x,yz] = [x,y][x,z][[z,x],y]  \label{eq-commutator_form_1}
\end{equation}
and
\begin{equation}
  [x^{-1},z]=[[x^{-1}, z],x][x,z]^{-1}, \hspace{1em} [x,y^{-1}] = [x,y]^{-1} [y,[y^{-1}, x]]. \label{eq-commutator_form_2}
\end{equation}
\item For subgroups $H$ and $K$ of $G$, we denote by $[H,K]$ the commutator subgroup of $G$
      generated by $[h, k]$ for all $h \in H$ and $k \in K$.
\item For elements $g_1, \ldots, g_k \in G$, a left-normed commutator
\[ [[ \cdots [[ g_{1},g_{2}],g_{3}], \cdots ], g_{k}] \in G \]
of weight $k$ is denoted by $[g_{1},g_{2}, \cdots, g_{k}]$.
\item For any $\Z$-module $A$, we denote the $\Q$-vector space $A \otimes_{\Z} \Q$ by the symbol obtained by attaching a subscript or superscript $\Q$ to $A$,
like $A_{\Q}$ or $A^{\Q}$. Similarly, for any $\Z$-linear map $f: A \rightarrow B$,
the induced $\Q$-linear map $A_{\Q} \rightarrow B_{\Q}$ is denoted by $f_{\Q}$ or $f^{\Q}$.
\end{itemize}

\section{Tangential and special derivation algebras of free Lie algebras}\label{S-TS-Der}

To begin with,  we review free Lie algebras and their derivation algebras.
For $n \geq 1$, let $F_n$ be the free group of rank $n$ with basis $x_1, \ldots , x_n$.
Let $\Gamma_n(1) \supset \Gamma_n(2) \supset \cdots$ be the lower central series of $F_n$ defined by
\[ \Gamma_n(1):= F_n, \hspace{1em} \Gamma_n(k) := [\Gamma_n(k-1),F_n], \hspace{1em} k \geq 2. \]
We denote by $\mathcal{L}_n(k) := \Gamma_n(k)/\Gamma_n(k+1)$ the $k$-th graded quotient of the lower central series of $F_n$,
and by $\mathcal{L}_n := {\bigoplus}_{k \geq 1} \mathcal{L}_n(k)$ the associated graded sum.
For each $k \geq 1$, we consider $\mathcal{L}_n(k)$ as an additive group. For example, for any $x, y, z \in F_n$
we have
\[ [xy,z]=[x,z]+[y,z], \hspace{1em} [x^{-1},y]=-[x,y] \]
in $\mathcal{L}_n(1)$ from (\ref{eq-commutator_form_1}) and (\ref{eq-commutator_form_2}).
The graded sum $\mathcal{L}_n$ has the Lie algebra structure with the Lie bracket induced from the commutator bracket of $F_n$.
It is known that the Lie algebra $\mathcal{L}_n$ is isomorphic to the free Lie algebra generated by $\mathcal{L}_n(1)$
by a classical work of Magnus. (For details, see \cite{MKS} and \cite{Reu} for example.)
Each of $\mathcal{L}_n(k)$ is a free abelian group, and its rank is given by
\begin{equation}\label{ex-witt}
 r_n(k) := \frac{1}{k} \sum_{d | k} \mu(d) n^{\frac{k}{d}}
\end{equation}
by a work of Witt \cite{Wit}, where $\mu$ is the M\"{o}bius function, and $d$ runs over all positive divisors of $k$.
Hall \cite{Hl1} constructed an explicit basis of $\mathcal{L}_n(k)$ with basic commutators of weight $k$.
For example, basic commutators of weight less than five are listed below.
\vspace{0.5em}
\begin{center}
{\renewcommand{\arraystretch}{1.3}
\begin{tabular}{|c|l|l|} \hline
  $k$  & basic commutators  &     \\ \hline
  $1$  & $x_i$ & $1 \leq i \leq n$                 \\ \hline
  $2$  & $[x_i,x_j]$        & $1 \leq j<i \leq n$  \\ \hline
  $3$  & $[x_i, x_j, x_l]$  & $i > j \leq l$    \\ \hline
  $4$  & $[x_i, x_j, x_l, x_m]$  & $i > j \leq l \leq m$  \\
       & $[[x_i, x_j], [x_l, x_m]]$  & $(i,j) > (l,m)$ \\ \hline
\end{tabular}}
\end{center}
\vspace{0.5em}
\noindent
Here $(i,j) > (l,m)$ means the lexicographic order on $\N \times \N$.
(For details, see \cite{Hl2} and \cite{Reu} for example.)

\vspace{0.5em}

For simplicity, we denote by $H$ the abelianization $F_n^{\mathrm{ab}}=H_1(F_n,\Z)$ of $F_n$.
The basis $x_1, \ldots, x_n$ of $F_n$ induces a basis of $H$ as a free abelian group. we also denote it by the same letters $x_1, \ldots, x_n$
by abuse of language.
By fixing this basis of $H$, we identify the automorphism group $\mathrm{Aut}\,H$ with the general linear group $\mathrm{GL}(n,\Z)$.
Then $\mathrm{GL}(n,\Z)$ naturally acts on each of $\mathcal{L}_n(k)$.
For example, as $\mathrm{GL}(n,\Z)$-modules we have $\mathcal{L}_n(1) = H$ and $\mathcal{L}_n(2) = \wedge^2 H$.
The universal enveloping algebra of $\mathcal{L}_n$ is the tensor algebra
\[ T(H):= \Z \oplus H \oplus H^{\otimes 2} \oplus \cdots \]
of $H$ over $\Z$. From Poincar\'{e}-Birkhoff-Witt's theorem, the natural map $\iota : \mathcal{L}_n \rightarrow T(H)$
is injective. We denote by $\iota_k : \mathcal{L}_n(k) \rightarrow H^{\otimes k}$ the degree $k$ part of $\iota$, and
consider $\mathcal{L}_n(k)$ as a $\mathrm{GL}(n,\Z)$-submodule of $H^{\otimes k}$ through $\iota_k$.
Similarly, we consider $\mathcal{L}_n^{\Q}(k)$ as a $\mathrm{GL}(n,\Q)$-submodule of $H_{\Q}^{\otimes k}$ through $\iota_k^{\Q}$.

\vspace{0.5em}

Next, we consider the derivation algebra of $\mathcal{L}_n$.
We denote it by
\[ \mathrm{Der}(\mathcal{L}_n) := \{ f : \mathcal{L}_n \xrightarrow{\Z-\mathrm{linear}} \mathcal{L}_n \,|\, f([x,y]) = [f(x),y]+ [x,f(y)],
    \,\,\, x, y \in \mathcal{L}_n \}. \]
The Lie algebra $\mathrm{Der}(\mathcal{L}_n)$ is a graded Lie algebra.
For $k \geq 0$, the degree $k$ part of $\mathrm{Der}(\mathcal{L}_n)$ is defined to be
\[ \mathrm{Der}(\mathcal{L}_n)(k) := \{ f \in \mathrm{Der}(\mathcal{L}_n) \,|\, f(x) \in \mathcal{L}_n(k+1), \,\,\, x \in H \}. \]
For $k \geq 1$, by considering the universality of the free Lie algebra, we identify $\mathrm{Der}(\mathcal{L}_n)(k)$ with
$\mathrm{Hom}_{\Z}(H,\mathcal{L}_n(k+1)) = H^* {\otimes}_{\Z} \mathcal{L}_n(k+1)$
as a $\mathrm{GL}(n,\Z)$-module, where $H^*$ means the $\Z$-linear dual group of $H$.
Let $x_1^*, \ldots, x_n^*$ be the dual basis of $x_1, \ldots, x_n \in H$. Then $\mathrm{Der}(\mathcal{L}_n)(k)$
is generated by $x_i^* \otimes [x_{j_1}, x_{j_2}, \ldots, x_{j_{k+1}}]$ for all $1 \leq i, j_1, \ldots, j_{k+1} \leq n$
as a $\Z$-module.
Let $\mathrm{Der}^+(\mathcal{L}_n)$ be the graded Lie subalgebra of $\mathrm{Der}(\mathcal{L}_n)$ with positive degrees.
Similarly, we define the graded Lie algebras $\mathrm{Der}(\mathcal{L}_{n}^{\Q})$ and $\mathrm{Der}^+(\mathcal{L}_{n}^{\Q})$ over $\Q$.
Then, we have $\mathrm{Der}(\mathcal{L}_{n}^{\Q}) = \mathrm{Der}(\mathcal{L}_n) \otimes_{\Z} \Q$, and
$\mathrm{Der}(\mathcal{L}_n^{\Q})(k)=H_{\Q}^* {\otimes}_{\Q} \mathcal{L}_n^{\Q}(k+1)$ for any $k \geq 1$.

\vspace{0.5em}

In this paper, we consider many calculations of Lie brackets of $\mathrm{Der}(\mathcal{L}_n)$ and its Lie subalgebras.
All calculations are based on the following formula.
\begin{equation}\label{eq-der-cal}
\begin{split}
   [x_i^* \otimes & [x_{i_1}, x_{i_2}, \ldots, x_{i_k}], x_j^* \otimes [x_{j_1}, x_{j_2}, \ldots, x_{j_l}]] \\
    & = \sum_{t=1}^k \delta_{j i_t} x_i^* \otimes [x_{i_1}, \ldots, x_{i_{t-1}}, [x_{j_1}, x_{j_2}, \ldots, x_{j_l}], x_{i_{t+1}}, \ldots , x_{i_k}] \\
    & \hspace{2em} - \sum_{s=1}^l \delta_{i j_s} x_j^* \otimes [x_{j_1}, \ldots, x_{j_{s-1}}, [x_{i_1}, x_{i_2}, \ldots, x_{i_k}], x_{j_{s+1}}, \ldots , x_{j_l}],
\end{split}
\end{equation}
where $\delta_{pq}$ means Kronecker's delta. It is perhaps easiest to understand the calculation if a simple example is exhibited. For example, we have
\[\begin{split}
  [x_1^* \otimes [x_1, x_2, x_3, x_2], x_2^* \otimes [x_3, x_1]]
    & = x_1^* \otimes [x_1, [x_3, x_1], x_3, x_2] + x_1^* \otimes [x_1, x_2, x_3, [x_3,x_1]] \\
    & \hspace{2em} - x_2^* \otimes [x_3, [x_1, x_2, x_3, x_2]]. 
  \end{split}\]

\begin{defn}
For any $k \geq 1$, let $\mathfrak{p}_n(k)$ be the $\Z$-submodule of $\mathrm{Der}(\mathcal{L}_n)(k)$ generated by elements
\[ x_i^* \otimes [x_{j_1}, \ldots, x_{j_k}, x_i] \]
for all $1 \leq i, j_1, \ldots, j_k \leq n$. Let $\mathfrak{p}_n$ be the graded sum $\bigoplus_{k \geq 1} \mathfrak{p}_n(k)$.
\end{defn}
\noindent
The Lie subalgebra $\mathfrak{p}_n$ is called the tangential derivation Lie algebra of $\mathcal{L}_n$.
Consider the symmetric group $\mathfrak{S}_n$ of degree $n$ as the subgroup of $\mathrm{GL}(n,\Z)$
consisting of permutations of the basis $x_1, \ldots, x_n$ of $H$. Then each of $\mathfrak{p}_n(k)$ is an
$\mathfrak{S}_n$-invariant $\Z$-submodule of $\mathrm{Der}(\mathcal{L}_n)(k)$.
We remark that
$x_i^* \otimes [x_j, x_i]$ for any $1 \leq i \neq j \leq n$ form a basis of $\mathfrak{p}_n(1)$, and for $k \geq 2$
$\mathfrak{p}_n(k)$ is isomorphic to $H \otimes_{\Z} \mathcal{L}_n(k)$ as an $\mathfrak{S}_n$-module by the correspondence
\[ x_i^* \otimes [x_{j_1}, \ldots, x_{j_k}, x_i] \mapsto x_i \otimes [x_{j_1}, \ldots, x_{j_k}]. \]
Hence, we have $\mathrm{rank}_{\Z}(\mathfrak{p}_n(1)) =n^2(n-1)/2$ and
$\mathrm{rank}_{\Z}(\mathfrak{p}_n(k)) = nr_n(k)$ for any $k \geq 2$.

\vspace{0.5em}

For any $k \geq 1$, an $n$-tuple $\alpha=(\alpha_1,\alpha_2, \ldots ,\alpha_n) \in \mathbb{Z}_{\ge 0}^n$
such that $\alpha_1+ \alpha_2+ \cdots +\alpha_n=k+1$ is called a composition of $k+1$, denoted by $\alpha \models k+1$.
For a composition $\alpha=(\alpha_1,\alpha_2, \ldots ,\alpha_n)$, if $\alpha_{i+1}=\cdots=\alpha_n=0$ then
we write $\alpha=(\alpha_1,\alpha_2, \ldots ,\alpha_i)$ for simplicity.
For any composition $\alpha \models k+1$,
let $p(\alpha)$ be the set of all permutations of $(1^{\alpha_1}, 2^{\alpha_2}, \cdots, n^{\alpha_n})$.
For example, if $\alpha =(2,1) \models 3$ then $p(\alpha)=\{ (1,1,2), (1,2,1), (2,1,1) \}$.
If a composition $\alpha=(\alpha_1,\alpha_2, \ldots ,\alpha_n)$ of $k+1$ satisfies $\alpha_1 \geq \alpha_2 \geq \cdots \geq \alpha_n$, then
$\alpha$ is called a partition of $k+1$, denoted by $\alpha \vdash k+1$.
Set
\[\begin{split}
   \mathfrak{p}_n(k,\alpha) := \langle x_{j_{k+1}}^* \otimes [x_{j_1},x_{j_2}, \ldots ,x_{j_{k+1}}]
      \,|\, (j_1, \ldots, j_{k+1}) \in p(\alpha) \rangle \subset \mathfrak{p}_n(k).
  \end{split}\]
Notice that each $\mathfrak{p}_n(k,\alpha)$ is not $\mathfrak{S}_n$-equivariant submodule of $\mathfrak{p}_n(k)$,
and that if $\alpha=(k+1)$ then $\mathfrak{p}_n(k,\alpha)=\{ 0 \}$ since $[x_1, x_1, \ldots, x_1]=0 \in \mathcal{L}_n(k+1)$.

\begin{lem}\label{L-TS-Der-1}
For any $k \geq 1$, we have
\[ \mathfrak{p}_n(k) = \bigoplus_{\alpha \vdash k+1} \mathfrak{S}_n \cdot \mathfrak{p}_n(k,\alpha), \]
where $\mathfrak{S}_n \cdot \mathfrak{p}_n(k,\alpha)$ means the $\mathfrak{S}_n$-orbit of $\mathfrak{p}_n(k,\alpha)$.\\
\end{lem}
\textit{Proof.}
It is clear that $\mathfrak{p}_n(k)$ is generated by $\mathfrak{S}_n \cdot \mathfrak{p}_n(k,\alpha)$
for all compositions $\alpha=(\alpha_1,\alpha_2, \ldots ,\alpha_n)$ of $k+1$.
In general, from the theory of commutator calculus given by Magnus and a theory of Hall basis,
a left normed commutator $[x_{i_1}, \ldots, x_{i_{k+1}}]$ can be written as a linear combination of basic commutators
among $x_{i_1}, \ldots, x_{i_{k+1}}$. In addition to this, all basic commutators are linearly independent.
(See Section in \cite{Hl2}.) This fact implies the required direct decomposition of $\mathfrak{p}_n(k)$.
$\square$

\vspace{0.5em}

For example, if $\alpha \models k+1$ such that $\alpha_1= \cdots = \alpha_i \neq 0$ and $\alpha_{i+1}= \cdots=\alpha_n=0$ for some
$1 \leq i \leq n$, then $\mathfrak{S}_n \cdot \mathfrak{p}_n(k,\alpha)$ is the direct sum of $\binom{n}{i}$ copies of $\mathfrak{p}_n(k,\alpha)$.
The above lemma is perhaps easy to understand if a few examples are given.
Consider the case of $k=1$. The compositions of $k+1=2$ are $(2)$ and $(1^2)$.
It suffices to consider the case for $(1^2)$ since $\mathfrak{p}_n(k,(2))=\{ 0 \}$.
If $\alpha=(1,1)$ then $(1^{\alpha_1}, 2^{\alpha_2}, \cdots, n^{\alpha_n})=(1,2)$, and we see
\begin{equation}\label{eq-p-1}
    \mathfrak{p}_n(1,\alpha) := \langle x_1^* \otimes [x_2,x_1], x_2^* \otimes [x_1, x_2] \rangle \subset \mathfrak{p}_n(1).
  \end{equation}
For any $1 \leq i<j \leq n$, set
\[ \mathfrak{p}_n(1,\alpha)_{i,j} := \langle x_i^* \otimes [x_j,x_i], x_j^* \otimes [x_i, x_j] \rangle \subset \mathfrak{p}_n(1). \]
Then we see
\[ \mathfrak{S}_n \cdot \mathfrak{p}_n(1,\alpha)= \bigoplus_{1 \leq i<j \leq n} \mathfrak{p}_n(1,\alpha)_{i,j} \cong \Z^{2 \binom{n}{2}}.  \]

\vspace{0.5em}

Next, we consider Lie algebras of special derivations of $\mathcal{L}_n$.
\begin{defn}
For any $k \geq 1$, let $\mathfrak{b}_n(k)$ be the $\mathfrak{S}_n$-invariant $\Z$-submodule
\[ \mathfrak{p}_n(k) \cap \Big{\{} \sum_{i=1}^n x_i^* \otimes C_i \in H^* \otimes_{\Z} \mathcal{L}_n(k+1) \,\Big{|}\, \sum_{i=1}^n C_i=0 \Big{\}} \]
of $\mathfrak{p}_n(k)$. Let $\mathfrak{b}_n$ be the graded sum $\bigoplus_{k \geq 1} \mathfrak{b}_n(k)$.
\end{defn}
\noindent
An element of $\mathfrak{b}_n$ is called a special derivation of $\mathcal{L}_n$.
It is also known that $\mathfrak{b}_n$ is a Lie subalgebra of $\mathfrak{p}_n$, and is called the special derivation Lie algebra of $\mathcal{L}_n$.
The rank of each $\mathfrak{b}_n(k)$ is given by Darn\'{e} \cite{Da3} in Section A.1 as
\begin{align}
 \mathrm{rank}_{\Z}(\mathfrak{b}_n(k)) 
    = \begin{cases}
        \binom{n}{2}, \hspace{1em} & k=1, \\
        nr_n(k) - r_n(k+1), & k \geq 2.
    \end{cases} 
\label{rank-bnk}
\end{align}

\vspace{0.5em}

For any $k \geq 1$ and any composition $\alpha=(\alpha_1,\alpha_2, \ldots ,\alpha_n) \models k+1$,
set $\mathfrak{b}_n(k,\alpha) = \mathfrak{p}_n(k,\alpha) \cap \mathfrak{b}_n(k)$.
From Lemma {\rmfamily \ref{L-TS-Der-1}}, we see
\begin{align}
 \mathfrak{b}_n(k) =
\bigoplus_{\alpha \vdash k+1} \mathfrak{S}_n \cdot \mathfrak{b}_n(k,\alpha), 
\label{decomp-bnk}
\end{align}
where $\mathfrak{S}_n \cdot \mathfrak{b}_n(k,\alpha)$ means the $\mathfrak{S}_n$-orbit of $\mathfrak{b}_n(k,\alpha)$.
From the above decomposition, we can calculate the dimension of $\mathfrak{b}_n^{\Q}(k)$ for any $k \geq 1$
by determining that of $\mathfrak{b}_n^{\Q}(k,\alpha)$ for each $\alpha \models k+1$ in principle.

\section{Johnson homomorphisms and trace maps}\label{S-John-trace}

Here we briefly recall the Andreadakis-Johnson filtrations and the Johnson homomorphisms of the automorphism groups of free groups, and their restrictions to the braid groups.
For each $k \geq 1$, the action of $\mathrm{Aut}\,F_n$ on the nilpotent quotient group $F_n/\Gamma_n(k+1)$ of $F_n$ induces the homomorphism
\[ \mathrm{Aut}\,F_n \rightarrow \mathrm{Aut}(F_n/\Gamma_n(k+1)). \]
We denote its kernel by $\mathcal{A}_n(k)$. Then the groups $\mathcal{A}_n(k)$ define a descending filtration
\[ \mathcal{A}_n(1) \supset \mathcal{A}_n(2) \supset \mathcal{A}_n(3) \supset \cdots \]
of $\mathrm{Aut}\,F_n$.
This filtration is called the Andreadakis-Johnson filtration of $\mathrm{Aut}\,F_n$, and the first term is called the IA-automorphism group
of $F_n$, and is also denoted by $\mathrm{IA}_n$.
Historically, the Andreadakis-Johnson filtration was originally introduced by Andreadakis \cite{And} in the 1960s.
Johnson studied such filtration for the mapping class groups of surfaces in the 1980s.
Andreadakis showed that
\begin{thm}[Andreadakis \cite{And}]\label{T-And} \quad
\begin{enumerate}
\item For any $k$, $l \geq 1$, $\sigma \in \mathcal{A}_n(k)$ and $x \in \Gamma_n(l)$, $x^{-1} x^{\sigma} \in \Gamma_n(k+l)$.
\item For any $k$ and $l \geq 1$, $[\mathcal{A}_n(k), \mathcal{A}_n(l)] \subset \mathcal{A}_n(k+l)$.
\end{enumerate}
\end{thm}
\noindent
For each $k \geq 1$, the group $\mathrm{Aut}\,F_n$ acts on $\mathcal{A}_n(k)$ by conjugation, and
this action naturally induces the actions of $\mathrm{GL}(n,\Z)=\mathrm{Aut}\,F_n/\mathrm{IA}_n$ on
the graded quotients $\mathrm{gr}^k (\mathcal{A}_n) := \mathcal{A}_n(k)/\mathcal{A}_n(k+1)$ by Part (2) of Theorem {\rmfamily \ref{T-And}}.

\vspace{0.5em}

In order to study the $\mathrm{GL}(n,\Z)$-module structure of ${\mathrm{gr}}^k (\mathcal{A}_n)$ for each $k \geq 1$,
we consider the Johnson homomorphisms.
For each $k \geq 1$, define the homomorphism
$\tilde{\tau}_k : \mathcal{A}_n(k) \rightarrow \mathrm{Hom}_{\Z}(H, {\mathcal{L}}_n(k+1))$ by
\[ \sigma \hspace{0.3em} \mapsto \hspace{0.3em} (x \mapsto x^{-1} x^{\sigma}), \hspace{1em} x \in H. \]
The kernel of $\tilde{\tau}_k$ is $\mathcal{A}_n(k+1)$. 
Hence, $\tilde{\tau}_k$ induces the injective homomorphism
\[ \tau_k : \mathrm{gr}^k (\mathcal{A}_n) \hookrightarrow \mathrm{Hom}_{\Z}(H, \mathcal{L}_n(k+1))
       = H^* \otimes_{\Z} \mathcal{L}_n(k+1). \]
The homomorphisms $\tilde{\tau}_k$ and ${\tau}_{k}$ are called the $k$-th Johnson homomorphisms of $\mathrm{Aut}\,F_n$.
Each of ${\tau}_{k}$ is $\mathrm{GL}(n,\Z)$-equivariant.
The graded sum $\mathrm{gr}(\mathcal{A}_n) := \bigoplus_{k \geq 1} \mathrm{gr}^k (\mathcal{A}_n)$ has the Lie algebra structure
with the Lie bracket induced from the commutator bracket on $\mathrm{IA}_n$.
Furthermore, the graded sum $\tau := \bigoplus_{k \geq 1} \tau_k : \mathrm{gr}(\mathcal{A}_n) \rightarrow \mathrm{Der}^+(\mathcal{L}_n)$
is a Lie algebra homomorphism.

\vspace{0.5em}

In the later sections, we consider restrictions of the Johnson homomorphisms to certain subgroups of $\mathrm{Aut}\,F_n$.
Let $\mathcal{S}$ be a subgroup of $\mathrm{Aut}\,F_n$. For any $k \geq 1$, set $\mathcal{S}(k):= \mathcal{A}_n(k) \cap \mathcal{S}$, and
$\mathrm{gr}^k(\mathcal{S}) := \mathcal{S}(k)/\mathcal{S}(k+1)$.
Then the restriction of $\tilde{\tau}_k$ to $\mathcal{S}(k)$ induces the injective homomorphism
$\tau_k^{\mathcal{S}} : \mathrm{gr}^k(\mathcal{S}) \rightarrow H^* \otimes_{\Z} \mathcal{L}_n(k+1)$.
We call $\tau_k^{\mathcal{S}}$ the $k$-th Johnson homomorphism of $\mathcal{S}$.
This homomorphism $\tau_k^{\mathcal{S}}$ is an $\mathcal{S}/\mathcal{S}(1)$-equivariant.
The graded sum $\mathrm{gr}(\mathcal{S}) = \bigoplus_{k \geq 1} \mathrm{gr}^k(\mathcal{S})$ is a Lie subalgebra of $\mathrm{gr}(\mathcal{A}_n)$, and
the graded sum $\tau := \bigoplus_{k \geq 1} \tau_k^{\mathcal{S}} : \mathrm{gr}(\mathcal{S}) \rightarrow \mathrm{Der}^+(\mathcal{L}_n)$ is a Lie algebra homomorphism.

\vspace{0.5em}

Here we remark on the cokernels of the Johnson homomorphisms.
For $k \geq 1$, let ${\varphi}^{k} : H^* {\otimes}_{\Z} H^{\otimes (k+1)} \rightarrow H^{\otimes k}$
be the contraction map defined by
\[ x_i^* \otimes x_{j_1} \otimes \cdots \otimes x_{j_{k+1}} \mapsto x_i^*(x_{j_1}) \, \cdot
    x_{j_2} \otimes \cdots \otimes \cdots \otimes x_{j_{k+1}}. \]
For the natural embedding ${\iota}_{k+1} : \mathcal{L}_n(k+1) \rightarrow H^{\otimes (k+1)}$,
we have the $\mathrm{GL}(n,\Z)$-equivariant homomorphism
\[ \Phi^k = {\varphi}^{k} \circ ({id}_{H^*} \otimes {\iota}_{k+1})
    : H^* {\otimes}_{\Z} \mathcal{L}_n(k+1) \rightarrow H^{\otimes k}. \]
For any $k \geq 1$,
let $\mathcal{C}_n(k)$ be the quotient module of $H^{\otimes k}$ by the action of the cyclic group $C_k$ of order $k$ on the components.
Namely,
\[ \mathcal{C}_n(k) = H^{\otimes k} \big{/} \langle a_1 \otimes a_2 \otimes \cdots \otimes a_k - a_2 \otimes a_3 \otimes \cdots \otimes a_k \otimes a_1
   \,|\, a_i \in H \rangle. \]
Let $\varpi^k : H^{\otimes k} \rightarrow \mathcal{C}_n(k)$ be the natural projection.
For any $1 \leq j_1, \ldots, j_k \leq n$, we write the image of $x_{j_1} \otimes \cdots \otimes x_{j_k} \in H^{\otimes k}$
by the map $\varpi^k$ as $x_{j_1} x_{j_2} \cdots x_{j_k} \in \mathcal{C}_n(k)$. According to this notation, we have
$x_{j_1} x_{j_2} \cdots x_{j_k}=x_{j_2} \cdots x_{j_k} x_{j_1}$ in $\mathcal{C}_n(k)$.
By \cite{Reu}, it is known that $\mathcal{C}_n(k)$ is a free abelian group, and
the necklaces among $x_1, \ldots, x_n$ of length $k$ form a basis of $\mathcal{C}_n(k)$.
From this, we have
\[ \mathrm{rank}_{\Z}(\mathcal{C}_n(k)) = \frac{1}{k}\sum_{d|k}\varphi(d)n^{\frac{k}{d}} \]
where $\varphi : \N \rightarrow \N$ is the Euler function, and $d$ runs over all positive divisors of $k$.

For any $k \geq 1$, set
$\mathrm{Tr}_k := \varpi^k \circ \Phi^k : H^* {\otimes}_{\Z} \mathcal{L}_n(k+1) \rightarrow \mathcal{C}_n(k)$.
We call $\mathrm{Tr}_k$ the trace map of degree $k$. This is a generalization of the Morita trace defined by Morita \cite{Mo1},
which is obtained by $\mathrm{Tr}_k$ by composing the natural projection $\mathcal{C}_n(k) \rightarrow S^kH$.
By combining independent results by Satoh \cite{S03, S13} and Darn\'{e} \cite{Dar}, for any $k \geq 2$ and $n \geq k+2$, the sequence
\[ 0 \rightarrow \mathrm{gr}^k (\mathcal{A}_n) \xrightarrow{\tau_k} H^* {\otimes}_{\Z} \mathcal{L}_n(k+1) \xrightarrow{\mathrm{Tr}_k} \mathcal{C}_n(k) \rightarrow 0 \]
is a $\mathrm{GL}(n,\Z)$-equivariant exact sequence. We remark that over $\Q$, the corresponding exact sequence was obtained by independent works by Satoh \cite{S03, S13} and
Massuyeau-Sakasai \cite{M-S}. In particular, we have
\begin{equation}
  \mathrm{Im}(\tau_k) = \mathrm{Ker}(\mathrm{Tr}_k) \label{eq-John-trace}
\end{equation}
for $n \geq k+2$. 

\section{McCool groups and pure braid groups}\label{S-M-PB}

Here we consider subgroups of the loop braid group $\mathrm{LB}_n$.
The intersection $\mathrm{LB}_n \cap \mathrm{IA}_n$ is
\[ \{ \sigma \in \mathrm{Aut}\,F_n \,|\,  x_i^{\sigma} = c_i x_i c_i^{-1} \hspace{0.5em} (c_i \in F_n) \}, \]
and called the basis-conjugating automorphisms of $F_n$, denoted by $\mathrm{P}\Sigma_n$.
McCool \cite{McC} obtained a finite presentation of $\mathrm{P}\Sigma_n$. Thus $\mathrm{P}\Sigma_n$ is also called the McCool group.
For any $1 \leq i \neq j \leq n$, let $K_{ij}$ be the basis-conjugating automorphism of $F_n$ defined by
\[ x_t \mapsto \begin{cases}
               {x_j}^{-1} x_i x_j & t=i, \\
               x_t                & t \neq i.
              \end{cases}\]
McCool obtained the following finite presentation of $\mathrm{P}\Sigma_n$.
\begin{thm}[McCool \cite{McC}]
The group $\mathrm{P}\Sigma_n$ is generated by $K_{ij}$ for $1 \leq i \neq j \leq n$ subject to relations: \\
\hspace{2em} {\bf{(P1)}}: $[K_{ij},K_{kj}]=1$ for $(i,j) < (k,j)$, \\
\hspace{2em} {\bf{(P2)}}: $[K_{ij}, K_{kl}]=1$ for $(i,j)< (k,l)$, \\
\hspace{2em} {\bf{(P3)}}: $[K_{ik}, K_{ij} K_{kj}]=1$ \\
where in each relation, the indices $i$, $j$, $k$, $l$ are distinct.
\end{thm}
\noindent
Let $\rho : \mathrm{Aut}\,F_n \rightarrow \mathrm{GL}(n,\Z)$ be the natural homomorphism induced from the abelianization of $F_n$.
The image of the restriction of $\rho$ to $\mathrm{LB}_n$ is the symmetric group $\mathfrak{S}_n$ of degree $n$.
Here we consider $\mathfrak{S}_n$ as a subgroup of $\mathrm{GL}(n,\Z)$, consisting of all permutations on the basis of $H$.
Thus we have the group extension
\[ 1 \rightarrow \mathrm{P}\Sigma_n \rightarrow \mathrm{LB}_n \xrightarrow{\rho} \mathfrak{S}_n \rightarrow 1. \]
For each $k \geq 1$, set $\mathrm{P}\Sigma_n(k) := \mathcal{A}_n(k) \cap \mathrm{P}\Sigma_n$, and $\mathrm{gr}^k(\mathrm{P}\Sigma_n):=\mathrm{P}\Sigma_n(k)/\mathrm{P}\Sigma_n(k+1)$.
The $k$-th Johnson homomorphism of $\mathrm{LB}_n$ is denoted by $\tau_k^{LB} : \mathrm{gr}^k(\mathrm{P}\Sigma_n) \rightarrow H^* \otimes_{\Z} \mathcal{L}_n(k+1)$.
In Proposition 3.5 in \cite{S18}, we showed that for each $k \geq 2$, the image of $\tau_k^{LB}$ is contained in $\mathfrak{p}_n(k)$.
From this viewpoint, in this paper we consider $\tau_k^{LB} : \mathrm{gr}^k(\mathrm{P}\Sigma_n) \rightarrow \mathfrak{p}_n(k)$.
Each of $\tau_k^{LB}$ is $\mathfrak{S}_n$-equivariant, and
the cokernel of $\tau_k^{LB}$, denoted by $\mathrm{Coker}(\tau_k^{LB})$, always means the quotient module $\mathfrak{p}_n(k)/\mathrm{Im}(\tau_k^{LB})$.
In Enomoto-Satoh \cite{ES3}, we studied $\mathrm{Coker}(\tau_k^{LB})$ by using the trace maps.
At the present stage, the range of $k$ where $\mathrm{Coker}(\tau_k^{LB})$ is completely determined for any $n \geq 3$ is $1 \leq k \leq 4$.

\vspace{0.5em}

Let $\mathrm{Tr}_k^{LB} : \mathfrak{p}_n(k) \rightarrow \mathcal{C}_n(k)$ be the restriction of the trace map
$\mathrm{Tr}_k : H^* \otimes_{\Z} \mathcal{L}_n(k+1) \rightarrow \mathcal{C}_n(k)$ to $\mathfrak{p}_n(k)$.
In Proposition 3.7 in \cite{S18}, we showed that 
$\langle x_1^k, \ldots, x_n^k \rangle \cap \mathrm{Im}(\mathrm{Tr}_k^{LB}) = \{ 0 \}$ for any $k \geq 2$.
Here $\langle x_1^k, \ldots, x_n^k \rangle$ denotes the $\Z$-linear span of $x_1^k, \ldots, x_n^k$ in $\mathcal{C}_n(k)$.
Set $\overline{\mathcal{C}}_n(k) := \mathcal{C}_n(k)/\langle x_1^k, \ldots, x_n^k \rangle$.
By compositing $\mathrm{Tr}_k^{LB}$ and
the natural projection $\mathcal{C}_n(k) \rightarrow \overline{\mathcal{C}}_n(k)$,
we obtain the map $\overline{\mathrm{Tr}}_k^{LB} : \mathfrak{p}_n(k) \rightarrow \overline{\mathcal{C}}_n(k)$.
We also call $\overline{\mathrm{Tr}}_k^{LB}$ the trace map.
We remark that $\overline{\mathrm{Tr}}_k^{LB}$ is not surjective in general. For example, it is the case for $k=4$.
In fact, $x_i^2x_j^2$ and $x_ix_jx_ix_j$ for $1 \leq i < j \leq n$ form a basis of $\mathrm{Coker}(\overline{\mathrm{Tr}}_4^{LB})$
as a $\Z$-module. (For details, see Lemma 7.7 in \cite{ES3}.)

\vspace{0.5em}

Next, we consider the plane braid group of $n$ strands.
Let $B_n$ be the mapping class group of $n$ points punctured disk $D^2 \setminus \{ p_1, \ldots, p_n \}$.
Namely, $B_n$ is the group of isotopy classes of orientation preserving homeomorphisms of $D^2 \setminus \{ p_1, \ldots, p_n \}$ which fix $\partial D^2$ pointwise.
By a classical work of Artin \cite{Art}, it is known that $B_n$ is isomorphic to the plane braid group of $n$ strands, and $B_n$
can be embedded into $\mathrm{Aut}\,F_n$. He \cite{Art} also gave the image of the embedding by
\[ \{ \sigma \in \mathrm{Aut}\,F_n \,|\, x_i^{\sigma} = c_i x_{\mu(i)} c_i^{-1} \hspace{0.5em} (c_i \in F_n, \,\,\, \mu \in \mathfrak{S}_n),
      \hspace{0.5em} (x_1x_2 \cdots x_n)^{\sigma}= x_1x_2 \cdots x_n \}. \]
Let $\sigma_1, \ldots, \sigma_{n-1}$ be Artin's generators of $B_n$. As an automorphism of $\mathrm{Aut}\,F_n$, $\sigma_i$ is defined by 
\[ x_t \mapsto \begin{cases}
               x_i^{-1} x_{i+1} x_i & t=i, \\
               x_i & t=i+1, \\
               x_t & t \neq i, i+1.
              \end{cases}\]

The intersection $B_n \cap \mathrm{IA}_n$ is called the pure braid group of $n$-strands, and is denoted by $P_n$. We have
\[ P_n = \{ \sigma \in \mathrm{Aut}\,F_n \,|\,  x_i^{\sigma} = c_i x_i c_i^{-1} \hspace{0.5em} (c_i \in F_n),
    \hspace{0.5em} (x_1x_2 \cdots x_n)^{\sigma}= x_1x_2 \cdots x_n \} \]
as a subgroup of $\mathrm{Aut}\,F_n$.
For any $1 \leq i<j \leq n$, set
\[ A_{i,j} := \sigma_{j-1} \sigma_{j-2} \cdots \sigma_{i+1} \sigma_i^2 \sigma_{i+1}^{-1} \cdots \sigma_{j-2}^{-1} \sigma_{j-1}^{-1} \in P_n. \]
Then Artin \cite{Art} showed that $P_n$ is generated by $A_{i,j}$ for all $1 \leq i<j \leq n$, and also obtained a finite presentation of $P_n$ with respect to these generators.
(See also Lemma 1.8.2 in \cite{Bir}.)
The image of the restriction of $\rho$ to $B_n$ is also the symmetric group $\mathfrak{S}_n$ in $\mathrm{GL}(n,\Z)$ as mentioned above.
Thus we have the following commutative diagram whose three rows are group extensions, and whose all vertical maps are inclusion maps.
\[\begin{CD}
1 @>>> \mathrm{IA}_n @>>> \mathrm{Aut}\,F_n @>{\rho}>> \mathrm{GL}(n,\Z) @>>> 1 \\
   @.     @AAA     @AAA      @AAA        @.  \\
1 @>>> \mathrm{P}\Sigma_n @>>> \mathrm{BP}_n @>{\rho|_{\mathrm{BP}_n}}>> \mathfrak{S}_n @>>> 1 \\
   @.     @AAA     @AAA      @AA{=}A        @.  \\
1 @>>> P_n @>>> B_n @>{\rho|_{B_n}}>> \mathfrak{S}_n @>>> 1 \\
 \end{CD}\]
For each $k \geq 1$, set $B_n(k) := \mathcal{A}_n(k) \cap B_n$, and $\mathrm{gr}^k(B_n):=B_n(k)/B_n(k+1)$.
Darn\'{e} \cite{Da2} showed that the filtration $B_n(1) \supset B_n(2) \supset \cdots$ coincides with the lower central series of the pure braid group $P_n$.
Hence from formulae obtained by Kohno \cite{Koh} and Falk-Randell \cite{F-R} independently, we can calculate the rank of each of the graded quotients $\mathrm{gr}^k(B_n)$.
For example, we have the following.
\begin{center}
{\renewcommand{\arraystretch}{1.3}
\begin{tabular}{|c|c|c|c|c|c|} \hline
  $k$  & $1$ & $2$ & $3$ & $4$ \\ \hline
  $\mathrm{rank}_{\Z}(\mathrm{gr}^k(B_n))$  & $\binom{n}{2}$ & $\binom{n}{3}$ & $2\binom{n+1}{4}$ & $\frac{3(2n-1)}{5}\binom{n+1}{4}$   \\ \hline
\end{tabular}}
\end{center}
\vspace{0.5em}
\noindent
The $k$-th Johnson homomorphism of $B_n$ is denoted by $\tau_k^B : \mathrm{gr}^k(B_n) \rightarrow H^* \otimes_{\Z} \mathcal{L}_n(k+1)$.
\begin{lem}\label{L-5-1}
For each $k \geq 2$, the image of $\tau_k^B$ is contained in $\mathfrak{b}_n(k)$.
\end{lem}
\textit{Proof.}
For any $\sigma \in B_n$, and any $1 \leq i \leq n$, set $x_i^{\sigma} := a_i x_i x_i$ for some $a_i \in F_n$.
Set $C_i := x_i^{-1} x_i^{\sigma} :=[x_i^{-1}, a_i^{-1}]$ for any $1 \leq i \leq n$.
From $(x_1 \cdots x_n)^{\sigma}=x_1 \cdots x_n$, we have
$(x_1C_1 \cdots x_nC_n)=x_1 \cdots x_n$.

If $\sigma \in B_n(k)$, we have $C_i \in \Gamma_n(k+1)$ for any $1 \leq i \leq n$.
Thus we see
\[ x_1 \cdots x_nC_1 \cdots C_n \equiv x_1C_1 \cdots x_nC_n = x_1 \cdots x_n \pmod{\Gamma_n(k+2)}, \]
and hence $C_1 \cdots C_n \equiv 1 \pmod{\Gamma_n(k+2)}$. From this, we obtain the required result.
$\square$

\vspace{0.5em}

From this viewpoint, in this paper we consider $\tau_k^B : \mathrm{gr}^k(B_n) \rightarrow \mathfrak{b}_n(k)$.
Each of $\tau_k^B$ is $\mathfrak{S}_n$-equivariant, and
the cokernel of $\tau_k^B$, denoted by $\mathrm{Coker}(\tau_k^B)$, always means the quotient module $\mathfrak{b}_n(k)/\mathrm{Im}(\tau_k^B)$.
Let $\mathrm{Tr}_k^B : \mathfrak{b}_n(k) \rightarrow \mathcal{C}_n(k)$ be the restriction of
$\mathrm{Tr}_k^{LB} : \mathfrak{p}_n(k) \rightarrow \mathcal{C}_n(k)$ to $\mathfrak{b}_n(k)$, and
$\overline{\mathrm{Tr}}_k^{B} : \mathfrak{b}_n(k) \rightarrow \overline{\mathcal{C}}_n(k)$
the composition of $\mathrm{Tr}_k^B$ and the natural projection $\mathcal{C}_n(k) \rightarrow \overline{\mathcal{C}}_n(k)$.
We emphasize the following open problem.
\begin{prob}
Determine $\mathrm{Im}(\overline{\mathrm{Tr}}_k^{LB})$ for any $k \geq 5$, and
$\mathrm{Im}(\overline{\mathrm{Tr}}_k^B)$ for any $k \geq 1$.
\end{prob}

\section{Basis of $\mathfrak{b}_n(k)$ for degrees $k \leq 4$}\label{S-strB}

Note that the rank of $\mathfrak{b}_n(k)$ is given in (\ref{rank-bnk}). 
In this section, we will show the following decompositions of $\mathfrak{b}_n(k)$ for any  $n \geq 2$ and $k \leq 4$ through (\ref{decomp-bnk});
\begin{align*}
\mathfrak{b}_n(1)&=\mathfrak{S}_n \cdot \mathfrak{b}_n(1,(1^2)),\\
\mathfrak{b}_n(2)&=\mathfrak{S}_n\cdot \mathfrak{b}_n(2,(1^3)),\\
\mathfrak{b}_n(3)&=\mathfrak{S}_n\cdot(\mathfrak{b}_n(3,(2^2))\oplus \mathfrak{b}_n(3,(2,1^2))\oplus \mathfrak{b}_n(3,(1^4))),\\
\mathfrak{b}_n(4)&=\mathfrak{S}_n\cdot(\mathfrak{b}_n(4,(3,1^2))\oplus \mathfrak{b}_n(4,(2^2,1))\oplus \mathfrak{b}_n(4,(2,1^3))\oplus \mathfrak{b}_n(4,(1^5))).
\end{align*}
Moreover, in order to investigate the abelianization of $\mathfrak{b}_n$, we give an explicit basis of each components. In this section, we will use the following special deriveations;
\begin{align*}
&x_i^*\otimes [x_j,x_i]+x_j^*\otimes [x_i,x_j] \in \mathfrak{b}_n(1), \\
&x_i^* \otimes [x_\ell,x_j,x_i]- x_j^* \otimes [x_\ell,x_i,x_j]+ x_\ell^* \otimes [x_j,x_i,x_\ell]\in \mathfrak{b}_n(2),  \\
&x_i^* \otimes [x_j,x_i,x_j,x_i] + x_j^* \otimes [x_i,x_j,x_i,x_j] \in \mathfrak{b}_n(3), \\
\bm{b}(i,j,l,p):=& x_i^* \otimes [x_l,x_j,x_p,x_i] + x_j^* \otimes [x_p,x_i,x_l,x_j] - x_l^* \otimes [x_p,x_i,x_j,x_l] \\
     & + x_p^* \otimes [x_j,x_l,x_i,x_p] \in \mathfrak{b}_n(3),\\
\bm{b}(i, j,l,p,q) :=& x_q^* \otimes [x_i,x_j,x_l,x_p,x_q] - x_p^* \otimes [x_i,x_j,x_l,x_q,x_p] - x_j^* \otimes [x_p,x_q,x_l,x_i,x_j] \\
   &+ x_i^* \otimes [x_p,x_q,x_l,x_j,x_i] - x_l^* \otimes [x_i,x_j,[x_p,x_q],x_l] \in \mathfrak{b}_n(4).
\end{align*}
We remark that Alekseev and Torossian \cite{AT1} gave a connection between trivalent trees and derivations of free Lie algebras, and that
\v{S}evera-Willwacher \cite{SW1} described $\mathfrak{b}_n$ as a space of undirected internally trivalent trees modulo the IHX relations.
(See Lemma 4 in \cite{SW1}.) We also refer a brief description of this given in Appendix A in \cite{Fel} by Felder.
From this point of view, the above five elements are described by the following internal trivalent trees, respectively.
(The third element is a half of the element corresponding to the third tree.)
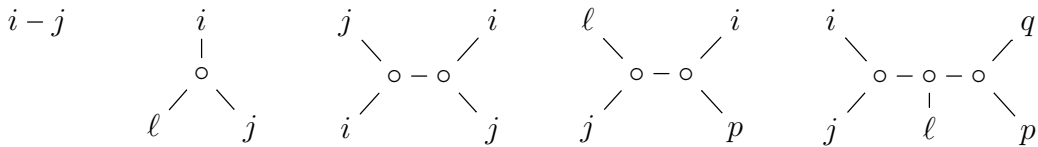
\[
\xymatrix@R=2mm@C=2mm{
i\ar@{-}[r] &j
} \quad \quad
\xymatrix@R=2mm@C=2mm{
& i\ar@{-}[d] & \cr
& \circ\ar@{-}[dr]\ar@{-}[dl] & \cr
\ell && j
} \quad \quad 
\xymatrix@R=2mm@C=2mm{
j\ar@{-}[rd] & & & i\ar@{-}[ld] \cr
& \circ\ar@{-}[r]&\circ & \cr
i\ar@{-}[ru] & & & j\ar@{-}[lu]
} \quad \quad
\xymatrix@R=2mm@C=2mm{
\ell\ar@{-}[dr]&&&i \cr
& \circ\ar@{-}[r] & \circ\ar@{-}[ru]\ar@{-}[rd] &\cr
j\ar@{-}[ru]&&&p
} \quad \quad
\xymatrix@R=2mm@C=2mm{
i\ar@{-}[dr]&&&&q \cr
& \circ\ar@{-}[r] & \circ\ar@{-}[r]\ar@{-}[d] & \circ\ar@{-}[ru]\ar@{-}[rd] &\cr
j\ar@{-}[ru]&&\ell&&p
}
\]

\subsection{The cases of $k=1$ and $2$}\label{Sss-1}
For $k=1$, the compositions of $k+1=2$ are $\alpha=(2)$ and $(1^2)$.
It suffices to consider the case for $(1^2)$ since $\mathfrak{b}_n(k,(2))=\{ 0 \}$.
By using (\ref{eq-p-1}), we can see
\[\begin{split}
    \mathfrak{b}_n(1, (1^2)) = \langle x_1^* \otimes [x_2,x_1] + x_2^* \otimes [x_1, x_2] \rangle \cong \Z.
  \end{split}\]
\noindent
For $k=2$, consider compositions $\alpha:=(2,1)$ and $(1^3)$ of $3$.
Since $\mathfrak{p}_n(2, (2,1))=\langle x_1^* \otimes [x_2,x_1,x_1] \rangle \cong \Z$, we see
$\mathfrak{b}_n(2, (2,1))=\{ 0 \}$. For $n \geq 3$, we have
\[ \mathfrak{p}_n(2, (1^3))=\langle x_1^* \otimes [x_3,x_2,x_1], x_2^* \otimes [x_3,x_1,x_2], x_3^* \otimes [x_2,x_1,x_3] \rangle \cong \Z^3. \]
Take any element $f = ax_1^* \otimes [x_3,x_2,x_1]+bx_2^* \otimes [x_3,x_1,x_2]+cx_3^* \otimes [x_2,x_1,x_3] \in \mathfrak{b}_n(2, (1^3))$.
By using the Jacobi identity, we have
\[\begin{split}
   0 &= a [x_3,x_2,x_1]+ b [x_3,x_1,x_2] + c[x_2,x_1,x_3] \\
     &= (a+b)[x_3,x_1,x_2] +(c-a)[x_2,x_1,x_3] \in \mathcal{L}_n(3).
  \end{split}\]
Since $[x_3,x_1,x_2]$ and $[x_2,x_1,x_3]$ are basic commutators of weight three, they are linearly independent in $\mathcal{L}_n(3)$,
and hence $a+b=c-a=0$.
This shows
\begin{equation}\label{eq-gen-b-2}
\begin{split}
    \mathfrak{b}_n(2,(1^3)) = \langle x_1^* \otimes [x_3,x_2,x_1]- x_2^* \otimes [x_3,x_1,x_2]+ x_3^* \otimes [x_2,x_1,x_3] \rangle \cong \Z.
\end{split}
\end{equation}

\subsection{The case of $k=3$}\label{Sss-3}
For $k=3$, consider compositions $\alpha=(3,1)$, $(2^2)$, $(2,1^2)$ and $(1^4)$ of $4$.
Since $\mathfrak{p}_n(3,(3,1))=\langle x_1^* \otimes [x_2,x_1,x_1,x_1] \rangle \cong \Z$, we see
$\mathfrak{b}_n(3,(3,1))=\{ 0 \}$. Since
\[ \mathfrak{p}_n(3,(2^2))=\langle x_1^* \otimes [x_2,x_1,x_2,x_1], x_2^* \otimes [x_1,x_2,x_1,x_2] \rangle \cong \Z^2, \]
we can see
\[ \mathfrak{b}_n(3,(2^2))=\langle x_1^* \otimes [x_2,x_1,x_2,x_1] + x_2^* \otimes [x_1,x_2,x_1,x_2] \rangle \cong \Z \]
since
\[ [x_2,x_1,x_2,x_1] + [x_1,x_2,x_1,x_2] = [x_2,x_1,[x_2,x_1]]=0 \in \mathcal{L}_n(4) \]
by the Jacobi identity.

\vspace{0.5em}

Assume $n \geq 3$.
For $\alpha=(2,1^2)$, we see that $\mathfrak{p}_n(3,(2,1^2))$ is generated by
\[ x_1^* \otimes [x_2,x_1,x_3,x_1], x_1^* \otimes [x_3,x_1,x_2,x_1], x_2^* \otimes [x_3,x_1,x_1,x_2], x_3^* \otimes [x_2,x_1,x_1,x_3]. \]
By using the Jacobi identity, we see
\[ x_1^* \otimes [x_2,x_1,x_3,x_1] + x_1^* \otimes [x_3,x_1,x_2,x_1] - x_2^* \otimes [x_3,x_1,x_1,x_2] - x_3^* \otimes [x_2,x_1,x_1,x_3]
   \in \mathfrak{b}_n(3,(2,1^2)) \]
since
\[\begin{split}
   & \underline{[x_2, x_1,x_3,x_1]} + [x_3,x_1,x_2,x_1] - [x_3,x_1,x_1,x_2] - \underline{[x_2,x_1,x_1,x_3]} \\
    & \hspace{0.5em} = \underline{[x_2, x_1, [x_3,x_1]]} + [x_3,x_1, [x_2,x_1]]=0.
  \end{split}\]
From this observation, we can see $\mathfrak{b}_n(3,(2,1^2)) = \langle \bm{b}(1,1,2,3) \rangle \cong \Z$
since $[x_3,x_1,x_2,x_1]$, $[x_3,x_1,x_1,x_2]$ and $[x_2,x_1,x_1,x_3]$ are linearly independent in $\mathcal{L}_n(4)$.

\vspace{0.5em}

For any $1 \leq i, j, l, p \leq n$, observe
\[\begin{split}
   & [x_l,x_j,x_p,x_i] + \underline{[x_p,x_i,x_l,x_j]} - \underline{[x_p,x_i,x_j,x_l]} - [x_l,x_i,x_j,x_p] + [x_j,x_i,x_l,x_p] \\
   &= [x_l,x_j,x_p,x_i] + \underline{[x_p,x_i,[x_l,x_j]]} - [x_l,x_j,x_i,x_p] \\
   &= [x_l,x_j,[x_p,x_i]] + [x_p,x_i,[x_l,x_j]]=0.
  \end{split}\]
\noindent
Assume $n \geq 4$.
For $\alpha=(1^4)$, we see that $\mathfrak{p}_n(3,(1^4))$ is generated by
\[\begin{split}
  & x_1^* \otimes [x_4,x_2,x_3,x_1], x_1^* \otimes [x_3,x_2,x_4,x_1], x_2^* \otimes [x_4,x_1,x_3,x_2], x_2^* \otimes [x_3,x_1,x_4,x_2], \\
  & x_3^* \otimes [x_4,x_1,x_2,x_3], x_3^* \otimes [x_2,x_1,x_4,x_3], x_4^* \otimes [x_3,x_1,x_2,x_4], x_4^* \otimes [x_2,x_1,x_3,x_4].
  \end{split}\]
We can see $\mathfrak{b}_n(3,(1^4)) = \langle \bm{b}(1,2,3,4), \bm{b}(1,2,4,3) \rangle \cong \Z^2$,
since the commutators appearing in the above except for
$[x_3,x_2,x_4,x_1]$ and $[x_4,x_2,x_3,x_1]$ are linearly independent in $\mathcal{L}_n(4)$.

\subsection{The case of $k=4$}\label{Sss-4}
For $k=4$, consider compositions $\alpha=(4,1)$, $(3,2)$, $(3,1^2)$,
$(2^2,1)$, $(2,1^3)$ and $(1^5)$ of $5$.
Since $\mathfrak{p}_n(4,(4,1))=\langle x_1^* \otimes [x_2,x_1,x_1,x_1,x_1] \rangle \cong \Z$, we see
$\mathfrak{b}_n(4,(4,1))=\{ 0 \}$. For $\alpha=(3,2)$, we have
\[ \mathfrak{p}_n(4,(3,2))=\langle x_1^* \otimes [x_2,x_1,x_1,x_2,x_1], x_2^* \otimes [x_2,x_1,x_1,x_1, x_2] \rangle \cong \Z^2. \]
Since
\[ [x_2,x_1,x_1,x_2,x_1] = [x_2,x_1,x_1,[x_2,x_1]] + [x_2,x_1,x_1,x_2,x_1], \]
and since $[x_2,x_1,x_1,[x_2,x_1]]$, $[x_2,x_1,x_1,x_2,x_1]$ and $[x_2,x_1,x_1,x_1, x_2]$ are basic commutators,
we can see $\mathfrak{b}_n(4,(3,2))=\{ 0 \}$.

For any $1 \leq i, j, l, p, q \leq n$, observe
\[\begin{split}
   & [x_i,x_j,x_l,x_p,x_q] -[x_i,x_j,x_l,x_q,x_p] - [x_p,x_q,x_l,x_i,x_j] + [x_p,x_q,x_l,x_j,x_i] \\
   & \hspace{0.5em} - [x_i,x_j,x_p,x_q,x_l] + [x_i,x_j,x_q,x_p,x_l] \\
  = & [x_i,x_j,x_l,[x_p,x_q]] - [x_p,x_q,x_l,[x_i,x_j]] - [x_i,x_j,[x_p,x_q],x_l] \\
  = & 0 \in \mathcal{L}_n(5).
  \end{split}\]
\noindent
It is easily seen that
\[ \bm{b}(i,j,l,p,q) = -\bm{b}(j,i,l,p,q) = -\bm{b}(i,j,l,q,p), \]
and
\begin{equation}\label{eq-b-r1}
\bm{b}(i,j,l,p,q)=\bm{b}(p,q,l,j,i).
\end{equation}
\noindent
Assume $n \geq 3$.
For $\alpha=(3,1^2)$, we see that $\mathfrak{p}_n(4,(3,1^2))$ is generated by
\[\begin{split}
  & x_1^* \otimes [x_3,x_1,x_1,x_2,x_1], x_1^* \otimes [x_2,x_1,x_1,x_3,x_1], x_1^* \otimes [x_3,x_1,[x_2,x_1],x_1], \\
  & x_2^* \otimes [x_3,x_1,x_1,x_1, x_2], x_3^* \otimes [x_2,x_1,x_1,x_1,x_3].
  \end{split}\]
Consider an element
\[\begin{split}
  \bm{b}(2, & 1,1,1,3) = x_3^* \otimes [x_2,x_1,x_1,x_1,x_3] - x_1^* \otimes [x_2,x_1,x_1,x_3,x_1] \\
   & \hspace{0.5em} - x_1^* \otimes [x_1,x_3,x_1,x_2,x_1] + x_2^* \otimes [x_1,x_3,x_1,x_1, x_2]- x_1^* \otimes [x_2,x_1,[x_1,x_3],x_1] 
  \end{split}\]
in $\mathfrak{b}_n(4,(3,1^2))$.
By an argument similar to as above, we can see $\mathfrak{b}_n(4,(3,1^2))=\langle \bm{b}(2,1,1,1,3) \rangle \cong \Z$.
We remark that $\bm{b}(3,1,1,1,2) =- \bm{b}(2,1,1,1,3)$.

\vspace{0.5em}

For $\alpha=(2^2,1)$, we see that $\mathfrak{p}_n(4,(2^2,1))$ is generated by
\[\begin{split}
  & x_1^* \otimes [x_2,x_1,x_2,x_3,x_1], x_1^* \otimes [x_3,x_1,x_2,x_2,x_1], x_1^* \otimes [x_3,x_2,[x_2,x_1],x_1], \\
  & x_2^* \otimes [x_3,x_1,x_1,x_2, x_2], x_2^* \otimes [x_2,x_1,x_1,x_3,x_2], x_2^* \otimes [x_3,x_1,[x_2,x_1],x_2], \\
  & x_3^* \otimes [x_2,x_1,x_1,x_2, x_3].
  \end{split}\]
Consider an element
\[\begin{split}
  \bm{b}(2, & 1,1,2,3) = x_3^* \otimes [x_2,x_1,x_1,x_2,x_3] - x_2^* \otimes [x_2,x_1,x_1,x_3,x_2] \\
   & \hspace{0.5em} - x_1^* \otimes [x_2,x_3,x_1,x_2,x_1] + x_2^* \otimes [x_2,x_3,x_1,x_1,x_2] - x_1^* \otimes [x_2,x_1,[x_2,x_3],x_1] 
  \end{split}\]
in $\mathfrak{b}_n(4,(2^2,1))$.
By an argument similar to as above, we can see $\mathfrak{b}_n(4,(2^2,1))=\langle \bm{b}(2,1,1,2,3) \rangle \cong \Z$.
We remark that an easy calculation shows that $\bm{b}(1,2,2,1,3)=-\bm{b}(2,1,1,2,3)$.

\vspace{0.5em}

Assume $n \geq 4$.
For $\alpha=(2,1^3)$, we see that $\mathfrak{p}_n(4,(2,1^3))$ is generated by
\[\begin{split}
  & x_1^* \otimes [x_4,x_1,x_2,x_3,x_1], x_1^* \otimes [x_3,x_1,x_2,x_4,x_1], x_1^* \otimes [x_2,x_1,x_3,x_4,x_1], \\
  & x_1^* \otimes [x_4,x_1,[x_3,x_2],x_1], x_1^* \otimes [x_4,x_3,[x_2,x_1],x_1], x_1^* \otimes [x_4,x_2,[x_3,x_1],x_1], \\
  & x_2^* \otimes [x_4,x_1,x_1,x_3, x_2], x_2^* \otimes [x_3,x_1,x_1,x_4,x_2], x_2^* \otimes [x_4,x_1,[x_3,x_1],x_2], \\
  & x_3^* \otimes [x_4,x_1,x_1,x_2, x_3], x_3^* \otimes [x_2,x_1,x_1,x_4, x_3], x_3^* \otimes [x_4,x_1,[x_2,x_1], x_3], \\
  & x_4^* \otimes [x_3,x_1,x_1,x_2, x_4], x_4^* \otimes [x_2,x_1,x_1,x_3, x_4], x_4^* \otimes [x_3,x_1,[x_2,x_1], x_4].
  \end{split}\]
Consider linearly independent elements
\[\begin{split}
  \bm{b}(1,2,1,3,4) & = -x_4^* \otimes [x_2,x_1,x_1,x_3,x_4] + X, \\
  \bm{b}(1,3,1,2,4) & = -x_4^* \otimes [x_3,x_1,x_1,x_2,x_4] + Y, \\
  \bm{b}(3,1,4,2,1) & = -x_4^* \otimes [x_3,x_1,[x_2,x_1],x_4] + Z
  \end{split}\]
in $\mathfrak{b}_n(4,(2,1^3))$, where $X, Y, Z$ mean linear combinations of elements type of $x_p^* \otimes [C, x_p]$ for $1 \leq p \leq 3$.
By a direct calculation, we can see that
\[\begin{split}
  & [x_4,x_1,x_2,x_3,x_1], [x_3,x_1,x_2,x_4,x_1], [x_2,x_1,x_3,x_4,x_1], \\
  & [x_4,x_1,[x_3,x_2],x_1], [x_4,x_3,[x_2,x_1],x_1], [x_4,x_2,[x_3,x_1],x_1], \\
  & [x_4,x_1,x_1,x_3, x_2], [x_3,x_1,x_1,x_4,x_2], [x_4,x_1,[x_3,x_1],x_2], \\
  & [x_4,x_1,x_1,x_2, x_3], [x_2,x_1,x_1,x_4, x_3], [x_4,x_1,[x_2,x_1], x_3], \\
  & [x_3,x_1,x_1,x_2, x_4], [x_2,x_1,x_1,x_3, x_4], [x_3,x_1,[x_2,x_1], x_4]
  \end{split}\]
are linearly independent in $\mathcal{L}_n(5)$, and hence
\[ \mathfrak{b}_n(4,(2,1^3)) = \langle \bm{b}(1,2,1,3,4), \bm{b}(1,3,1,2,4), \bm{b}(3,1,4,2,1) \rangle \cong \Z^3. \]
We remark that $\bm{b}(1,4,1,2,3) = \bm{b}(1,2,1,3,4) - \bm{b}(1,3,1,2,4)- 2 \bm{b}(3,1,4,2,1)$.
Hence,
\[ \mathfrak{b}_n^{\Q}(4,(2,1^3)) = \langle \bm{b}(1,2,1,3,4), \bm{b}(1,3,1,2,4), \bm{b}(1,4,1,2,3) \rangle_{\Q} \cong \Q^3. \]

\vspace{0.5em}

Assume $n \geq 5$.
Consider the case of $\alpha=(1^5)$. For each $1 \leq i \leq 5$,
let $1 \leq i_1 < i_2 < i_3 < i_4 \leq 5$ be the integers such that
$\{ 1, 2, 3, 4, 5 \} \setminus \{ i \} =\{ i_1, i_2, i_3, i_4 \}$.
Then $\mathfrak{p}_n(4,(1^5))$ is generated by
\begin{equation}\label{eq-muri}
\begin{split}
   & x_i^* \otimes [x_{i_4},x_{i_1},x_{i_2},x_{i_3}, x_i], x_i^* \otimes [x_{i_3},x_{i_1},x_{i_2},x_{i_4}, x_i], x_i^* \otimes [x_{i_2},x_{i_1},x_{i_3},x_{i_4}, x_i], \\
   & x_i^* \otimes [x_{i_4},x_{i_3},[x_{i_2},x_{i_1}], x_i], x_i^* \otimes [x_{i_4},x_{i_2},[x_{i_3},x_{i_1}], x_i], x_i^* \otimes [x_{i_4},x_{i_1},[x_{i_2},x_{i_3}], x_i]
\end{split}
\end{equation}
for all $1 \leq i \leq 5$ and $\mathfrak{p}_n(4,(1^5)) \cong \Z^{30}$.
Consider six linearly independent elements
{\small
\[\begin{split}
  \bm{b}(4,1,2,3,5) &= x_5^* \otimes [x_4,x_1,x_2,x_3,x_5] + X_1, \hspace{0.5em} \bm{b}(4,3,5,2,1) = -x_5^* \otimes [x_4,x_3,[x_2,x_1],x_5] + X_2, \\
  \bm{b}(3,1,2,4,5) &= x_5^* \otimes [x_3,x_1,x_2,x_4,x_5] + X_3, \hspace{0.5em} \bm{b}(4,2,5,3,1) = -x_5^* \otimes [x_4,x_2,[x_3,x_1],x_5] + X_4, \\
  \bm{b}(2,1,3,4,5) &= x_5^* \otimes [x_2,x_1,x_3,x_4,x_5] + X_5, \hspace{0.5em} \bm{b}(4,1,5,3,2) = -x_5^* \otimes [x_4,x_1,[x_3,x_2],x_5] + X_6,  
  \end{split}\]
}
\noindent
in $\mathfrak{b}_n(5,(1^5))$, where $X_m$ for each $1 \leq m \leq 6$ means a linear combination of elements type of $x_p^* \otimes [C, x_p]$ for $1 \leq p \leq 4$.
Then by a computer calculation,
we can see that elements in $(\ref{eq-muri})$ except for the elements for $i=5$ are linearly independent.
(We can also check it by a brute force hand calculation by using basic commutators of weight $5$.)
Thus $\mathfrak{b}_n(4,(1^5)) \cong \Z^{6}$. 

\section{Trace maps on $\mathfrak{b}_n(k)$ and the Johnson cokernels}\label{Ss-tarceB}

In this section, for $k \leq 4$ we describe the cokernels $\mathrm{Coker}(\tau_k^B)=\mathfrak{b}_n(k)/\mathrm{Im}(\tau_k^B)$ 
by using trace maps on $\mathfrak{b}_n(k)$.

\vspace{0.5em}

First we consider the case where $k=1$ and $2$. In these cases, the Johnson images are completely determined by Darn\'{e} \cite{Da3}.
So we briefly describe the facts along the lines of our arguments.
For $k=1$, by using the Artin generators $A_{i,j}$ of the pure braid group $P_n=B_n(1)$, we see that
\[ \tau_1^B(A_{i,j})=x_i^* \otimes [x_j, x_i] + x_j^* \otimes [x_i, x_j] \]
for any $1 \leq i<j \leq n$.
Hence from the result in Subsection {\rmfamily \ref{Sss-1}}, we obtain
$\mathrm{Im}(\tau_1^B)=\mathfrak{S}_n \cdot \mathfrak{b}_n(1,(1^2)) = \mathfrak{b}_n(1)$.
By using Artin's finite presentation of $P_n$, we see that
the abelianization of $P_n$ is a free abelian group of rank $\binom{n}{2}$ with freely generated by the coset classes of
the Artin generators. Thus, we also see that $\tau_1^B$ gives the abelianization of $P_n$.
Similarly, for $k=2$, we have
\[\begin{split}
 \tau_2^B([A_{ij},A_{il}]) = x_i^* \otimes [x_j, x_l,x_i] + x_j^* \otimes [x_l, x_i,x_j] + x_l^* \otimes [x_i, x_j,x_l]
 \end{split}\]
for any $1 \leq i, j, l \leq n$ such that $i< j, l$.
Hence  from the result in Section {\rmfamily \ref{Sss-1}}, we obtain
$\mathrm{Im}(\tau_2^B)=\mathfrak{S}_n \cdot \mathfrak{b}_n(2, (1^3)) = \mathfrak{b}_n(2)$.

\vspace{0.5em}

For $k \geq 3$, in order to study $\mathrm{Coker}(\tau_k^B)$, we consider Morita's trace maps.
For any $k \geq 1$, let $S^k H$ be the symmetric tensor product of $H$ of degree $k$.
The composition of $\mathrm{Tr}_k : H^* \otimes_{\Z} \mathcal{L}_n(k+1) \rightarrow \mathcal{C}_n(k)$ and
the natural projection $\mathcal{C}_n(k) \rightarrow S^k H$
is called the Morita trace map, and denoted by $\mathrm{MT}_k : H^* \otimes_{\Z} \mathcal{L}_n(k+1) \rightarrow S^k H$.
The restriction of $\mathrm{MT}_k$ to $\mathfrak{b}_n(k)$ is denoted by $\mathrm{MT}_k^B : \mathfrak{b}_n(k) \rightarrow S^k H$.
We show that the image of $\mathrm{MT}_k^B$ is non-trivial for any odd $k \geq 1$.
For any $1 \leq i_1, \ldots, i_k \leq n$, set $\bm{e}_{i_1, i_2, \ldots, i_k} := x_{i_1} x_{i_2} \cdots x_{i_k} \in S^k H$.
Then 
\[ \{ \bm{e}_{i_1, i_2, \ldots, i_k} \in S^k H \,|\, 1 \leq i_1 \leq \cdots \leq i_k \leq n \} \]
forms a basis of $S^k H$ as a free abelian group.
With respect to the map $\Phi^k$ discussed in Section {\rmfamily \ref{S-John-trace}},
from Lemma 3.1 in \cite{S03}, we have the following.
\begin{lem}\label{L-S03-1}
For any $1 \leq i, i_1, \ldots, i_{k+1} \leq n$,
\[\begin{split}
   \Phi^k(x_i^* \otimes & [x_{i_1}, x_{i_2}, \ldots, x_{i_{k+1}}]) \\
     & = \delta_{i, i_1} x_{i_2} \otimes \cdots \otimes x_{i_{k+1}}
         - \sum_{l=2}^{k+1} \delta_{i, i_l} [x_{i_1}, \ldots, x_{i_{l-1}}] \otimes x_{i_{l+1}} \otimes \cdots \otimes x_{i_{k+1}}.
  \end{split} \]
In particular, we have
\[ \mathrm{MT}_k(x_i^* \otimes [x_i, x_{i_1}, \ldots, x_{i_{k}}]) = \bm{e}_{i_1, i_2, \ldots, i_k}. \]
\end{lem}

For any $k \geq 3$, let $\mathcal{S}_k$ be the set of $\bm{e}_{i_1, i_2, \ldots, i_k} \in S^k H$ such that $1 \leq i_1 \leq \cdots \leq i_k \leq n$ and
there are at least three distinct indices among $i_1, \ldots, i_k$.
\begin{pro}\label{P-strB-1}
For any $n \geq 3$ and any odd $k \geq 3$, we have $2\mathcal{S}_k \subset \mathrm{Im}(\mathrm{MT}_k^B)$. 
\end{pro}
\textit{Proof.}
For any indeces $1 \leq j, i_1, i_2, \ldots, i_{k-1} \leq n$ such that $j, i_1, i_2$ are distinct, consider
\[ [x_{i_1}, \ldots, x_{i_{k-1}}, [x_j, x_{i_1}]] + [x_j, x_{i_1}, [x_{i_1}, \ldots, x_{i_{k-1}}]] = 0 \in \mathcal{L}_n(k+1), \]
and
\[\begin{split}
  x_{i_1}^* & \otimes [x_{i_1}, x_{i_2}, \ldots, x_{i_{k-1}}, x_j, x_{i_1}] - x_j^* \otimes [x_{i_1}, x_{i_2}, \ldots, x_{i_{k-1}}, x_{i_1}, x_j] \\
  & + x_{i_{k-1}}^* \otimes [x_j, x_{i_1}, [x_{i_1}, \ldots, x_{i_{k-2}}], x_{i_{k-1}}] \\
%  & - x_{i_{k-2}}^* \otimes [x_j, x_{i_1}, x_{i_{k-1}}, [x_{i_1}, \ldots, x_{i_{k-3}}], x_{i_{k-2}}] \\
  & + \cdots \\
  & +(-1)^{k-4} x_{i_3}^* \otimes [x_j, x_{i_1}, x_{i_{k-1}}, \ldots, x_{i_4}, [x_{i_1}, x_{i_2}], x_{i_3}] \\
  & +(-1)^{k-3} x_{i_2}^* \otimes [x_j, x_{i_1}, x_{i_{k-1}}, \ldots, x_{i_3}, x_{i_1}, x_{i_2}] \\
  & +(-1)^{k-2} x_{i_1}^* \otimes [x_j, x_{i_1}, x_{i_{k-1}}, \ldots, x_{i_3}, x_{i_2}, x_{i_1}]
 \end{split}\]
in $\mathfrak{b}_n(k)$.
The image of the above element by $\mathrm{MT}_k^B$ is $2 \bm{e}_{j, i_1, i_2, \ldots, i_{k-1}}$ since $k$ is odd.
This shows the proposition.
$\square$

We continue to investigate the image of the Morita trace for a general odd degree $k \geq 3$.
Set
\[ \mathcal{T}_k' := \{ \bm{e}_{i,i,j, i_3, \ldots, i_{k-1}} + \bm{e}_{j,j,i, i_3, \ldots, i_{k-1}} \,|\, 1 \leq i, j, i_3, \ldots, i_{k-1} \leq n, \hspace{0.5em} i \neq j \}. \]
\begin{pro}\label{P-strB-1-2}
For any $n \geq 3$ and $k \geq 3$, we have $\mathcal{T}_k' \subset \mathrm{Im}(\mathrm{MT}_k^B)$.
\end{pro}
\textit{Proof.}
For any indices $1 \leq i, j, i_3, \ldots, i_{k-1} \leq n$ such that $i \neq j$, consider a Jacobi identity
\[\begin{split}
  [x_i, x_j, x_{i_3}, \ldots, & x_{i_{k-1}}, x_i, x_j] + [x_i, x_j, [x_i, x_j, x_{i_3}, \ldots, x_{i_{k-1}}]] \\
   &+ [x_j, [x_i, x_j, x_{i_3}, \ldots, x_{i_{k-1}}], x_i] = 0 \in \mathcal{L}_n(k+1).
  \end{split}\]
Since
\[\begin{split}
   [x_i, x_j, & [x_i, x_j, x_{i_3}, \ldots, x_{i_{k-1}}]] \\
    & = \sum_{l=1}^{k-3} (-1)^{l-1} [x_i, x_j, x_{i_{k-1}}, \ldots, x_{i_{k-l+1}}, [x_i, x_j, x_{i_3}, \ldots, x_{i_{k-l-1}}], x_{i_{k-l}}],
 \end{split}\]
we see
\begin{equation}\label{eq-1109}
\begin{split}
  x_j^* \otimes & [x_i, x_j, x_{i_3}, \ldots, x_{i_{k-1}}, x_i, x_j] \\
   & +\sum_{l=1}^{k-3} (-1)^{l-1} x_{i_{k-l}}^* \otimes [x_i, x_j, x_{i_{k-1}}, \ldots, x_{i_{k-l+1}}, [x_i, x_j, x_{i_3}, \ldots, x_{i_{k-l-1}}], x_{i_{k-l}}] \\
   & + x_i^* \otimes [x_j, [x_i, x_j, x_{i_3}, \ldots, x_{i_{k-1}}], x_i]
  \end{split}
\end{equation}
in $\mathfrak{b}_n(k)$.
The image of the above element by $\mathrm{MT}_k^B$ is 
$-(\bm{e}_{i,i,j, i_3, \ldots, i_{k-1}} + \bm{e}_{j,j,i, i_3, \ldots, i_{k-1}})$.
This shows the proposition.
$\square$

In the notation $\bm{e}_{i_1, i_2, \ldots, i_k}$, if $i_1= \cdots = i_l=i$ and $i_{l+1}= \cdots = i_k=j$ for some $1 \leq l \leq k$ then we write it like
$\bm{e}_{i^l, j^{k-l}}$. Set
\[ \mathcal{T}_k := \{ \bm{e}_{i^l,j^{k-l}} + \bm{e}_{i^{l+1},j^{k-l-1}} \,|\, 1 \leq i < j \leq n, \hspace{0.5em} 1 \leq l \leq k-2 \} \subset \mathcal{T}_k'. \]
Since $2\mathcal{S}_k$ and $\mathcal{T}_k$ are linearly independent in $S^k H$ for odd $k \geq 3$, we have
\[ \mathrm{rank}_{\Z}(\mathrm{Im}(\mathrm{MT}_{k}^B)) \geq \binom{n+k-1}{k} - \dfrac{n(n+1)}{2}. \]
In order to show the equality holds in the above equation, we show the following.

\begin{pro}\label{P-strB-1-3}
Let $n \geq 3$ and $k \geq 3$ be integers. For any $1 \leq i < j \leq n$ and any $a \in \Z \setminus \{ 0 \}$,
we have $a \bm{e}_{i,j^{k-1}} \notin \mathrm{Im}(\mathrm{MT}_k^B)$.
\end{pro}
\textit{Proof.}
For some $1 \leq i < j \leq n$ and some $a \in \Z \setminus \{ 0 \}$, assume $a \bm{e}_{i,j^{k-1}} \in \mathrm{Im}(\mathrm{MT}_k^B)$ and
$\mathrm{MT}_k^B(X) = a \bm{e}_{i,j^{k-1}}$ for $X \in \mathfrak{b}_n(k) \subset \mathfrak{p}_n(k)$.
Since elements $x_p^* \otimes [x_{i_1}, x_{i_2}, \ldots, x_{i_k}, x_{p}]$ for $1 \leq p, i_1, \ldots, i_k \leq n$ generate
$\mathfrak{p}_n(k)$, and since
\[ \mathrm{MT}_k(x_p^* \otimes [x_{i_1}, x_{i_2}, \ldots, x_{i_k}, x_{p}]) = (\delta_{p,i_1}-\delta_{p, i_2}) \bm{e}_{i_1, i_2, \ldots, i_k} \]
from Lemma {\rmfamily \ref{L-S03-1}}, we may assume that $X$ is a linear combination of elements
$x_p^* \otimes [x_{i_1}, x_{i_2}, \ldots, x_{i_k}, x_{p}]$ such that $i_1, \ldots, i_k$ is eqaul to one $i$ and $k-1$ times $j$s, and $p=i$ or $j$.
Namely, since $i_1$ and $i_2$ should be different, $X$ is a linear combination of $k$ elements
\[\begin{split}
   & x_i^* \otimes [x_i, x_j, \ldots, x_j, x_i], x_j^* \otimes [x_i, x_j, \ldots, x_j], \\
   & x_j^* \otimes [x_i, x_j, \ldots, x_j, x_i, x_j], \ldots, x_j^* \otimes [x_i, x_j, x_i, x_j, \ldots, x_j].
  \end{split}\]
On the other hand, by reducing $X$ with the element (\ref{eq-1109}) of $\mathfrak{b}_n(k)$ in Proposition {\rmfamily \ref{P-strB-1-2}},
we may assume $X$ is a linear combination of $k-1$ elements
\[ x_j^* \otimes [x_i, x_j, \ldots, x_j], x_j^* \otimes [x_i, x_j, \ldots, x_j, x_i, x_j], \ldots, x_j^* \otimes [x_i, x_j, x_i, x_j, \ldots, x_j]. \]
Recall $X \in \mathfrak{b}_n(k)$, and consider when a linear combination of
\[ [x_i, x_j, \ldots, x_j], [x_i, x_j, \ldots, x_j, x_i, x_j], \ldots, [x_i, x_j, x_i, x_j, \ldots, x_j] \]
is equal to zero in $\mathcal{L}_n(k+1)$. Since $[x_i, x_j, \ldots, x_j]$ is the only element among the above whose
components have one $x_i$ and $k-1$ $x_j$s, by considering the Hall basis of $\mathcal{L}_n(k+1)$, we see that
the coefficient of $[x_i, x_j, \ldots, x_j]$ in such a linear combination is equal to zero.
Hence $X$ is a linear combination of
\[ x_j^* \otimes [x_i, x_j, \ldots, x_j, x_i, x_j], \ldots, x_j^* \otimes [x_i, x_j, x_i, x_j, \ldots, x_j]. \]
This is a contradiction to $\mathrm{MT}_k^B(X)=a \bm{e}_{i,j^{k-1}}$.
$\square$

From Proposition {\rmfamily \ref{P-strB-1-3}}, we see that the intersection of
$\mathrm{Im}(\mathrm{MT}_k^B)$ and the $\Z$-submodule of $S^k H$ generated by
\[ \{ \bm{e}_{i^k} \,|\, 1 \leq i \leq n \} \cup \{ \bm{e}_{i,j^{k-1}} \,|\, 1 \leq i<j \leq n \} \]
is trivial. Thus we obtain the following.
\begin{thm}\label{T-im-mtr}
For any $n \geq 3$ and any odd $k \geq 3$, we have
\[ \mathrm{rank}_{\Z}(\mathrm{Im}(\mathrm{MT}_{k}^B)) = \binom{n+k-1}{k} - \dfrac{n(n+1)}{2}. \]
\end{thm}

For $k=3$, we show that the equality holds. More precisely, we show the following.
\begin{pro}\label{T-John-1}
For any $n \geq 3$,
\[ 0 \rightarrow \mathrm{gr}^3(B_n) \xrightarrow{\tau_3^B} \mathfrak{b}_n(3) \xrightarrow{\mathrm{MT}_3^B} \mathrm{Im}(\mathrm{MT}_3^B) \rightarrow 0 \]
is an $\mathfrak{S}_n$-equivariant exact sequence.
\end{pro}
\textit{Proof.}
First we determine the image $\mathrm{Im}(\mathrm{MT}_3^B)$ of the Morita trace.
By Proposition {\rmfamily \ref{P-strB-1}}, we see that 
$\mathrm{MT}_3^B(\mathfrak{S}_n \cdot \mathfrak{b}_n(3,(2,1^2)))=2 \langle \mathcal{S}_3 \rangle = 2 \langle \bm{e}_{i,j,l} \,|\, 1 \leq i<j<l \leq n \rangle \cong \Z^{\binom{n}{3}}$.
Moreover, we can easily see
\[ \mathrm{MT}_3^B(\mathfrak{S}_n \cdot \mathfrak{b}_n(3,(2^2))) = \langle \bm{e}_{i,i,j} + \bm{e}_{j,j,i} \,|\,
    1 \leq i < j \leq n \rangle \cong \Z^{\binom{n}{2}}, \]
and $\mathrm{MT}_3^B(\mathfrak{S}_n \cdot \mathfrak{b}_n(3,(1^4)))=\{ 0 \}$. Hence we have
\[ \mathrm{Im}(\mathrm{MT}_3^B)=\mathrm{MT}_3^B(\mathfrak{S}_n \cdot (\mathfrak{b}_n(3,(2,1^2)) \oplus \mathfrak{b}_n(3,(2^2))))
    \cong \Z^{\binom{n+1}{3}}. \]

Next, we show $\mathrm{Im}(\tau_3^B) = \mathrm{Ker}(\mathrm{MT}_3^B)$. 
If this equation is shown, we obtain the required exact sequence since $\mathrm{Im}(\mathrm{MT}_3^B)$ is a free abelian group,.
It suffices to show $\mathrm{Im}(\tau_3^B) \supset \mathrm{Ker}(\mathrm{MT}_3^B)$
since $\mathrm{Im}(\tau_3^B) \subset \mathrm{Ker}(\mathrm{MT}_3^B)$.
Since 
\[ \mathrm{Im}(\mathrm{MT}_3^B)= \mathrm{MT}_3^B(\mathfrak{S}_n \cdot \mathfrak{b}_n(3,(2^2)))
     \oplus \mathrm{MT}_3^B(\mathfrak{S}_n \cdot \mathfrak{b}_n(3,(2,1^2))), \]
we see
\[\begin{split}
   \mathrm{Ker}(\mathrm{MT}_3^B) &= \Big{(} (\mathfrak{S}_n \cdot \mathfrak{b}_n(3,(2^2))) \cap \mathrm{Ker}(\mathrm{MT}_3^B) \Big{)}
    \oplus \Big{(}(\mathfrak{S}_n \cdot \mathfrak{b}_n(3,(2,1^2))) \cap \mathrm{Ker}(\mathrm{MT}_3^B) \Big{)} \\
    & \hspace{2em} \oplus \mathfrak{S}_n \cdot \mathfrak{b}_n(3,(1^4)). 
  \end{split}\]
From a straightforward calculation, we see
\[ \tau_3^B([A_{1,4}, [A_{1,3},A_{1,2}]])=\bm{b}(1,2,3,4), \hspace{0.5em} \tau_3^B([A_{1,3}, [A_{1,4},A_{1,2}]])=\bm{b}(1,2,4,3), \]
and hence $\mathfrak{S}_n \cdot \mathfrak{b}_n(3,(1^4)) \subset \mathrm{Im}(\tau_3^B)$.
Notice that $\mathfrak{S}_n \cdot \mathfrak{b}_n(3,(2^2))$ is generated by $-x_i^* \otimes [x_i,x_j,x_j,x_i] + x_j^* \otimes [x_i,x_j,x_j,x_i]$
for all $1 \leq j<i \leq n$ as a $\Z$-module. It turns out that $(\mathfrak{S}_n \cdot \mathfrak{b}_n(3,(2^2))) \cap \mathrm{Ker}(\mathrm{MT}_3^B)$
is trivial from
\[ \mathrm{MT}_3^B(-x_i^* \otimes [x_i,x_j,x_j,x_i] + x_j^* \otimes [x_i,x_j,x_j,x_i]) = -(\bm{e}_{i,i,j}+\bm{e}_{j,j,i}) \]
for any $1 \leq j<i \leq n$.
Notice that $\mathfrak{S}_n \cdot \mathfrak{b}_n(3,(2,1^2))$ is generated by $\bm{b}(i,i,j,l)$
for all distinct $1 \leq i, j, l \leq n$ as a $\Z$-module. From the observation
\[\begin{split}
   \bm{b}(i,i,l,j)&=-\bm{b}(i,i,j,l)), \\
   \tau_3^B([A_{i,j},[A_{i,j},A_{i,l}]])&= \bm{b}(i,i,j,l) + \bm{b}(j,j,i,l)
  \end{split}\]
for any distinct $i, j$ and $l$, we see that any element $X$ of $\mathfrak{S}_n \cdot \mathfrak{b}_n(3,(2,1^2))$ can be written as
\[ X \equiv \sum_{i<j<l} a_{i,j,l} \bm{b}(i,i,j,l) \pmod{\mathrm{Im}(\tau_3^B)} \]
for some $a_{i,j,l} \in \Z$. If $X \in \mathrm{Ker}(\mathrm{MT}_3^B)$, then we see $X \equiv 0 \pmod{\mathrm{Im}(\tau_3^B)}$ by using
$\mathrm{MT}_3^B(\bm{b}(i,i,j,l))=-2\bm{e}_{i,j,k}$.
Therefore we obtain $\mathrm{Im}(\tau_3^B) \supset \mathrm{Ker}(\mathrm{MT}_3^B)$. $\square$

\vspace{0.5em}

For $k=4$, we consider another trace maps.
For any $k \geq 1$, let $\psi^k : H^* {\otimes}_{\Z} H^{\otimes (k+1)} \rightarrow H^{\otimes k}$
be the contraction map with respect to the first and the third components, defined by
\[ x_i^* \otimes x_{j_1} \otimes \cdots \otimes x_{j_{k+1}} \mapsto x_i^*(x_{j_2}) \, \cdot
    x_{j_1} \otimes x_{j_3} \otimes \cdots \otimes \cdots \otimes x_{j_{k+1}}, \]
and set $\Psi^k := \psi^k \circ (\mathrm{id}_{H^*} \otimes \iota_{k+1}) : H^* {\otimes}_{\Z} \mathcal{L}_n(k+1) \rightarrow H^{\otimes k}$.
\begin{lem}\label{L-strB-4}
For any $1 \leq i, i_1, \ldots, i_{k+1} \leq n$,
\[\begin{split}
   \Psi^k(x_i^* & \otimes [x_{i_1}, x_{i_2}, \ldots, x_{i_{k+1}}]) \\
     & = \delta_{i, i_2} x_{i_1} \otimes x_{i_3} \otimes \cdots \otimes x_{i_{k+1}} \\
     & \hspace{1em} - \sum_{l=2}^{k+1} x_{i_l} \otimes \Big{(} \delta_{i, i_1} x_{i_2} \cdots \otimes x_{i_{l-1}}
         \otimes x_{i_{l+1}} \otimes \cdots \otimes x_{i_{k+1}} \\
     & \hspace{1em} - \sum_{m=2}^{l-1} \delta_{i,i_m} [x_{i_1}, \ldots, x_{i_{m-1}}]
         \otimes x_{i_{m+1}}\otimes \cdots \otimes x_{i_{l-1}} \otimes x_{i_{l+1}} \otimes \cdots \otimes x_{i_{k+1}} \Big{)}.
  \end{split} \]
\end{lem}
\textit{Proof.}
Since
\[\begin{split}
   \Psi^k(x_i^* & \otimes [x_{i_1}, x_{i_2}, \ldots, x_{i_{k+1}}]) \\
     & = \delta_{i, i_2} x_{i_1} \otimes x_{i_3} \otimes \cdots \otimes x_{i_{k+1}} \\
     & \hspace{1em} - \sum_{l=2}^{k+1} x_{i_l} \otimes \Phi^{k-1}(x_i^* \otimes [x_{i_1}, \ldots, x_{i_{l-1}}] \otimes x_{i_{l+1}} \otimes \cdots \otimes x_{i_{k+1}}),
 \end{split} \]
we obtain the required result by Lemma {\rmfamily \ref{L-S03-1}}.
$\square$

\vspace{0.5em}

Set
\begin{equation} \label{eq-def-tr_4}
 \mathrm{Tr}_{[2,1^{k-2}]} := (id_{H} \otimes f_{[1^{k-1}]}) \circ \Psi^k : H^* {\otimes}_{\Z} \mathcal{L}_n(k+1)
                  \rightarrow H {\otimes}_{\Z} \wedge^{k-1} H,
\end{equation}
where $f_{[1^k]} : H^{\otimes k} \rightarrow \wedge^k H$ is the natural projection defined by
\[ f_{[1^k]}(x_{i_1} \otimes \cdots \otimes x_{i_k}) = x_{i_1} \wedge \cdots \wedge x_{i_k}. \]
Let $I$ be the $\mathrm{GL}(n,\Z)$-submodule of $H {\otimes}_{\Z} \wedge^{k-1} H$ defined by
\[\begin{split}
   I &:= \langle x \otimes z_1 \wedge \cdots \wedge z_{k-2} \wedge y + y \otimes z_1 \wedge \cdots \wedge z_{k-2} \wedge x 
            \,\,|\,\, x, y, z_t \in H \rangle.
  \end{split}\]
We see $\mathrm{rank}_{\Z} I =(k-1)\binom{n+1}{k}$.
In our previous paper \cite{S03}, we showed that $\mathrm{Im}(\tau_k) \subset \mathrm{Ker}(\mathrm{Tr}_{[2,1^{k-1}]})$ and
$\mathrm{Im}\,(\mathrm{Tr}_{[2,1^{k-1}]}^{\Q}) = I_{\Q}$ if $k \geq 4$ is even and $n \geq k+1$.
By using Theorem 9 in \cite{EKS}, we see that $\mathrm{Ker}(\mathrm{Tr}_k) \subset \mathrm{Ker}(\mathrm{Tr}_{[2,1^{k-1}]})$
for $n \geq k+2$.
By using Pieri's formula, we have
$H_{\Q} {\otimes}_{\Z} \wedge^{k-1} H_{\Q} \cong H_{\Q}^{[2,1^{k-2}]} \oplus \wedge^{k} H_{\Q}$.
If $k$ is even, since $I_{\Q}$ is contained in the kernel of the natural map
$H_{\Q} {\otimes}_{\Z} \wedge^{k-1} H_{\Q} \rightarrow \wedge^k H_{\Q}$ defined by
$x \otimes y_1 \wedge \cdots \wedge y_{k-1} \mapsto x \wedge y_1 \wedge \cdots \wedge y_{k-1}$,
we see that $I_{\Q}$ is the $\mathrm{GL}(n,\Q)$ irreducible polynomial module associated to the Young diagram $[2,1^{k-2}]$.

\vspace{0.5em}

Let $J \subset I$ be the submodule of $I$ generated by
\[ w(i,j;i_1, \ldots, i_{k-2}) := x_i \otimes x_{i_1} \wedge \cdots \wedge x_{i_{k-2}} \wedge x_j + x_j \otimes x_{i_1} \wedge \cdots \wedge x_{i_{k-2}} \wedge x_i \]
for all distinct $1 \leq i, j, i_1, \ldots, i_{k-2} \leq n$.
Then we show the following.
\begin{pro}\label{P-strB-3}
For any $k \equiv 0 \pmod{4}$ and $n \geq k$, we have $2^{\frac{k-2}{2}}J \subset \mathrm{Im}(\mathrm{Tr}_{[2,1^{k-2}]}|_{\mathfrak{b}_n(k)})$.
\end{pro}
\textit{Proof.}
For any distinct $1 \leq i, j, i_1, \ldots, i_{k-2} \leq n$, observe
\[\begin{split}
   [x_{i_1}, x_{i_2}, & x_i, [x_{i_3}, x_{i_4}], \ldots, [x_{i_{k-3}}, x_{i_{k-2}}], [x_j,x_i]] \\
   & + [x_j,x_i, [x_{i_1}, x_{i_2}, x_i, [x_{i_3}, x_{i_4}], \ldots, [x_{i_{k-3}}, x_{i_{k-2}}]]] =0 \in \mathcal{L}_n(k+1),
  \end{split}\]
and consider
\[\begin{split}
   x_i^* & \otimes [x_{i_1}, x_{i_2}, x_i, [x_{i_3}, x_{i_4}], \ldots, [x_{i_{k-3}}, x_{i_{k-2}}], x_j, x_i] \\
    & - x_j^* \otimes [x_{i_1}, x_{i_2}, x_i, [x_{i_3}, x_{i_4}], \ldots, [x_{i_{k-3}}, x_{i_{k-2}}], x_i, x_j] \\
    & + x_{i_{k-2}}^* \otimes [x_j,x_i, [x_{i_1}, x_{i_2}, x_i, [x_{i_3}, x_{i_4}], \ldots, [x_{i_{k-5}}, x_{i_{k-4}}]], x_{i_{k-3}}, x_{i_{k-2}}] \\
    & - x_{i_{k-3}}^* \otimes [x_j,x_i, [x_{i_1}, x_{i_2}, x_i, [x_{i_3}, x_{i_4}], \ldots, [x_{i_{k-5}}, x_{i_{k-4}}]], x_{i_{k-2}}, x_{i_{k-3}}] \\
    & \cdots \\
    & +(-1)^{\frac{k-6}{2}} x_{i_4}^* \otimes [x_j,x_i, [x_{i_{k-3}}, x_{i_{k-2}}], \ldots,  [x_{i_5}, x_{i_6}], [x_{i_1}, x_{i_2}, x_i], x_{i_3}, x_{i_4}] \\
    & +(-1)^{\frac{k-4}{2}} x_{i_3}^* \otimes [x_j,x_i, [x_{i_{k-3}}, x_{i_{k-2}}], \ldots,  [x_{i_5}, x_{i_6}], [x_{i_1}, x_{i_2}, x_i], x_{i_4}, x_{i_3}] \\
    & +(-1)^{\frac{k-4}{2}} x_i^* \otimes [x_j,x_i, [x_{i_{k-3}}, x_{i_{k-2}}], \ldots,  [x_{i_3}, x_{i_4}], [x_{i_1}, x_{i_2}], x_i] \\
    & +(-1)^{\frac{k-2}{2}} x_{i_2}^* \otimes [x_j,x_i, [x_{i_{k-3}}, x_{i_{k-2}}], \ldots,  [x_{i_3}, x_{i_4}], x_i, x_{i_1}, x_{i_2}] \\
    & +(-1)^{\frac{k}{2}} x_{i_1}^* \otimes [x_j,x_i, [x_{i_{k-3}}, x_{i_{k-2}}], \ldots,  [x_{i_3}, x_{i_4}], x_i, x_{i_2}, x_{i_1}] \\
  \end{split}\]
in $\mathfrak{b}_n(k)$.
By using Lemma {\rmfamily \ref{L-strB-4}}, we can calculate the image of the above element by $\mathrm{Tr}_{[2,1^{k-2}]}$ as follows:
\[\begin{split}
  x_i & \otimes (2 x_{i_1} \wedge x_{i_2}) \wedge \cdots \wedge (2x_{i_{k-3}} \wedge x_{i_{k-2}}) \wedge x_j \\
      & + x_j \otimes (2 x_{i_1} \wedge x_{i_2}) \wedge \cdots \wedge (2x_{i_{k-3}} \wedge x_{i_{k-2}}) \wedge x_i \\
      & +(-1)^{\frac{k-4}{2}} x_i \otimes x_j \wedge (2 x_{i_{k-3}} \wedge x_{i_{k-2}}) \wedge \cdots \wedge (2x_{i_1} \wedge x_{i_2}) \\
      & +(-1)^{\frac{k-4}{2}} x_j \otimes (2 x_{i_{k-3}} \wedge x_{i_{k-2}}) \wedge \cdots \wedge (2x_{i_1} \wedge x_{i_2}) \wedge x_i \\
      & = 2^{\frac{k-2}{2}} w(i,j;i_1, \ldots, i_{k-2}) 
 \end{split}\]
since $k \equiv 0 \pmod{4}$. This shows the proposition.
$\square$

\begin{pro}\label{T-John-2}
For any $n \geq 4$,
\[ 0 \rightarrow \mathrm{gr}_{\Q}^4(B_n) \xrightarrow{\tau_{4,\Q}^B} \mathfrak{b}_n^{\Q}(4) \xrightarrow{\mathrm{Tr}_{[2,1^2]}^{\Q}}
    I_{\Q} \rightarrow 0 \]
is an $\mathfrak{S}_n$-equivariant exact sequence.
\end{pro}
\textit{Proof.}
%First we determine $\mathrm{Im}(\mathrm{Tr}_{[2,1^2]}^{\Q}|_{\mathfrak{b}_n^{\Q}(4)})$. \mathrm{Im}(\mathrm{Tr}_{[2,1^2]}|_{\mathfrak{b}_n(4)})
By Proposition {\rmfamily \ref{P-strB-3}}, we have that $2^{\frac{k-2}{2}}J \subset \mathrm{Im}(\mathrm{Tr}_{[2,1^2]}|_{\mathfrak{b}_n(4)})$.
In addition to this, we have
\[ \mathrm{Tr}_{[2,1^2]}(\bm{b}(2,1,1,1,3)) = 4(x_1 \otimes x_2 \wedge x_3 \wedge x_1), \]
and hence
\[ \mathrm{Im}(\mathrm{Tr}_{[2,1^2]}^{\Q}|_{\mathfrak{b}_n^{\Q}(4)})=J_{\Q} \oplus
    \langle x_i \otimes x_j \wedge x_l \wedge x_i \,|\, 1 \leq i, j, l \leq n \rangle_{\Q}
  = I_{\Q}. \]
%In particular, $\mathrm{Im}(\mathrm{Tr}_{[2,1^2]}|_{\mathfrak{b}_n(4)}) \cong \Z^{3\binom{n+1}{4}}$.
Thus we see
\[ \mathrm{dim}_{\Q}(\mathrm{Im}(\tau_{4,\Q}^B)) + \mathrm{dim}_{\Q}(\mathrm{Im}(\mathrm{Tr}_{[2,1^2]}^{\Q}|_{\mathfrak{b}_n^{\Q}(4)}))
 = \mathrm{dim}_{\Q}(\mathfrak{b}_n^{\Q}(4)), \]
and obtain the required result. $\square$

%We remark that by tensoring with $\Q$, we have an $\mathfrak{S}_n$-equivariant exact sequence
%\[ 0 \rightarrow \mathrm{gr}^4(B_n^{\Q}) \xrightarrow{\tau_{4,\Q}^B} \mathfrak{b}_n^{\Q}(4) \xrightarrow{\mathrm{Tr}_{[2,1^2]}^{\Q}}
%    I_{\Q} \rightarrow 0. \]

\section{On the abelianization of $\mathfrak{b}_n$}\label{S-abel}

Here we consider the abelianization of the special derivation algebra $\mathfrak{b}_n$ as a Lie algebra.
The purpose of this section is to detect non-trivial linearly independent elements in $H_1(\mathfrak{b}_n,\Z)$ with
the Morita traces and others.

\subsection{Infinite generation of $H_1(\mathfrak{b}_n,\Z)$}\label{Ss-ab-inf}

In \cite{Mo3}, Morita constructed a surjective Lie algebra homomorphism
\[ \Theta := \mathrm{id}_1 \oplus \bigoplus_{k \geq 2} \mathrm{MT}_k : H_1(\mathrm{Der}^+(\mathcal{L}_n),\Z) \rightarrow (H^* \otimes_{\Z} \wedge^2 H) \oplus
   \bigoplus_{k \geq 2} S^k H \]
where $\mathrm{id}_1$ is the identity map on the degree one part
$\mathrm{Der}^+(\mathcal{L}_n)(1)=H^* \otimes_{\Z} \wedge^2 H$ of $\mathrm{Der}^+(\mathcal{L}_n)$,
and the target is understood to be an abelian Lie algebra.

\vspace{0.5em}

Since $\mathfrak{b}_n$ is a graded Lie algebra, its abelianization $H_1(\mathfrak{b}_n,\Z)$ naturally has a graded $\Z$-module structure.
We denote the degree $k$ part of $H_1(\mathfrak{b}_n,\Z)$ by $H_1(\mathfrak{b}_n,\Z)_k$ for each $k \geq 1$.
As a corollary to Theorem {\rmfamily \ref{T-im-mtr}}, by observing the restriction of $\Theta$ to $\mathfrak{b}_n$, we have the following.
\begin{thm}\label{T-strB-2}
For any $n \geq 3$ and odd $k \geq 3$, we have
\[ \mathrm{rank}_{\Z}(H_1(\mathfrak{b}_n,\Z)_k) \geq \binom{n+k-1}{k} - \dfrac{n(n+1)}{2}. \]
In particular, $H_1(\mathfrak{b}_n,\Z)$ is not finitely generated as a $\Z$-module.
\end{thm}
In Theorem {\rmfamily \ref{T-ninj-MT}}, we will show that $\Theta|_{\mathfrak{b}_n}$ is not injective.
Namely, the equality in Theorem {\rmfamily \ref{T-strB-2}} does not hold in general.

\begin{rem}\label{R-kuno}
Yusuke Kuno communicated to us the non-finite generation of $H_1(\mathfrak{b}_2,\Z)$ can be obtained with Soul\'{e} elements defined by Ihara \cite{Iha} as follows.
In a series of his study of the absolute Galois group over $\Q$, he considered a certain infinite sequence of solvable Galois
extensions over $\Q$ and studied the graded Lie algebra structure of the graded sum $\mathcal{G}:=\bigoplus_{k \geq 1} \mathcal{G}^{(k)}$ consisting of the Galois groups of
the successive extensions of the sequence.
In particular, he \cite{Iha} constructed a sequence of elements $\sigma_{2m+1} \in \mathcal{G}^{(2m+1)}$ for each $m \geq 1$ and called them Soul\'{e} elements.
In Equation (15) in \cite{AT2}, for any $m \geq 1$
Alekseev and Torossian \cite{AT2} described the images of $\sigma_{2m+1}$ by some Lie algebra homomorphism $\nu : \mathcal{G} \rightarrow \mathfrak{b}_2^{\Q}$.
More precisely, in Proposition 4.5 in \cite{AT2}, they showed that
$\mathrm{MT}_{2m+1}^{\Q}(\nu(\sigma_{2m+1}))=(x_1+x_2)^{2m+1}-x_1^{2m+1}-x_2^{2m+1}$.
This shows that the non-finite generation of
$H_1(\mathfrak{b}_2^{\Q},\Q)$ as a Lie algebra. Hence, by considering the natural map $\mathfrak{b}_2 \rightarrow \mathfrak{b}_2^{\Q}$, we see the non-finite generation of
$H_1(\mathfrak{b}_2,\Z)$ as a Lie algebra.
Futhermore, for any $n \geq 3$ and distinct $1 \leq i, j \leq n$,
by considering $\binom{n}{2}$ embeddings $\mathfrak{b}_2^{\Q} \rightarrow \mathfrak{b}_n^{\Q}$ induced by the correspondence
$x_1 \mapsto x_i$ and $x_2 \mapsto x_j$ respectively,
we can see that $\mathrm{dim}_{\Q}(\mathrm{Im}(\mathrm{MT}_{2m+1}^{\Q})) \geq \binom{n}{2}$. This shows that non-finite generation of
$H_1(\mathfrak{b}_n^{\Q},\Q)$, and of $H_1(\mathfrak{b}_n,\Z)$.
\end{rem}

\subsection{The structure of $H_1(\mathfrak{b}_n,\Z)_k$ for $k \leq 4$}\label{Ss-ab-low}

Here we consider the structure of $H_1(\mathfrak{b}_n,\Z)_k$ for $k \leq 4$.
Clearly, $H_1(\mathfrak{b}_n,\Z)_1 = \mathfrak{b}_n(1) \cong \Z^{\binom{n}{2}}$.
In this subsection, we show that $H_1(\mathfrak{b}_n,\Z)_k$ for $2 \leq k \leq 3$ and
$H_1(\mathfrak{b}_n,\Q)_4$ are trivial except for the components which can be detected by the Morita trace map.

\vspace{0.5em}

For $k=2$, since $\mathfrak{b}_n(2)= \mathfrak{S}_n \cdot \mathfrak{b}_n(2,(1^3))$ is generated by
$\tau_2^B([A_{1,3},A_{1,2}])$ as an $\mathfrak{S}_n$-module from (\ref{eq-gen-b-2}), we see $H_1(\mathfrak{b}_n,\Z)_2=0$.
For $k=3$, since $\mathrm{gr}(B_n):=\bigoplus_{k \geq 1} \mathrm{gr}^k(B_n)$ is generated by
the degree 1 part as a Lie algebra, it is immediately follows from Theorem {\rmfamily \ref{T-John-1}} that
$H_1(\mathfrak{b}_n,\Z)_3 \cong \mathrm{Im}(\mathrm{MT}_3^B) \cong \Z^{\binom{n+1}{3}}$.

\vspace{0.5em}

For $k=4$, we show that all generators $\bm{b}(i,j,l,p,q)$ of $\mathfrak{b}_n(4)$
belong to $[\mathfrak{b}_n(3),\mathfrak{b}_n(1)]$.
Notice that $\mathrm{Im}(\tau_4^B) \subset [\mathfrak{b}_n(3),\mathfrak{b}_n(1)]$, and
for any distinct $1 \leq i, j, l, p \leq n$,
\[ \bm{b}(i,j,l,p) \in \mathrm{Ker}(\mathrm{MT}_3^B) = \mathrm{Im}(\tau_3^B) \]
from Lemma {\rmfamily \ref{L-S03-1}}.
\begin{lem}
$(1)$\,\,For any distinct $1 \leq i, j, l, p, q \leq n$, we have $\bm{b}(i,j,l,p,q) \in \mathrm{Im}(\tau_4^B)$. \\
$(2)$\,\,For any distinct $1 \leq i, j, l \leq n$, we have $\bm{b}(i,j,j,i,l) \in [\mathfrak{b}_n(3),\mathfrak{b}_n(1)]$.
\end{lem}
\textit{Proof.}
This lemma is obtained by straightforward calculations given by
\[\begin{split}
   \bm{b}(i,j,l,p,q) &= [\bm{b}(l,q,p,j), \tau_1^B(A_{p,q})], \\
   \bm{b}(i,j,j,i,l) &= [x_j^* \otimes [x_i, x_j, x_i, x_j] + x_i^* \otimes [x_j, x_i, x_j, x_i], \tau_1^B(A_{i,l})].
  \end{split}\]
$\square$

\vspace{0.5em}

Let us consider $\bm{b}(i,j,j,j,l), \bm{b}(i,j,j,p,q)$ for any distinct $1 \leq i, j, l, p \leq n$.
\begin{lem}\label{L-John-7}
For distinct $1 \leq i,j,l \leq n$, we have\\
$(1)$\,\,$\bm{b}(i,j,j,i,l) \in [\mathfrak{b}_n(3),\mathfrak{b}_n(1)]$. \\
$(2)$\,\,$\bm{b}(i,j,j,j,l) - \bm{b}(i,j,l,l,j) \in [\mathfrak{b}_n(3),\mathfrak{b}_n(1)]$. \\
$(3)$\,\,$\bm{b}(i,j,j,j,l) \in [\mathfrak{b}_n(3),\mathfrak{b}_n(1)]$.
\end{lem}
\textit{Proof.} Parts (1) and (2) are follow from
\[\begin{split}
   \bm{b}(i,j,j,i,l) &= [x_i^* \otimes [x_j,x_i,x_j,x_i] + x_j^* \otimes [x_i,x_j,x_i,x_j], \tau_1^B(A_{i,l})], \\
   \bm{b}(i,j,j,j,l) &-\bm{b}(i,j,l,l,j) = [\tau_1^B(A_{j,l}), \bm{b}(j,j,i,l)].
  \end{split}\]
From (\ref{eq-b-r1}), we see $\bm{b}(j,i,l,l,j)=\bm{b}(i,l,l,j,i)=-\bm{b}(i,l,l,i,j)$. Thus
by Parts (1) and (2), we obtain Part (3). $\square$

\vspace{0.5em}

\begin{lem}\label{L-John-9}
For distinct $1 \leq i, j, l, p \leq n$, $\bm{b}(i,j,j,p,q) \in [\mathfrak{b}_n(3),\mathfrak{b}_n(1)]$.
\end{lem}
\textit{Proof.} 
This lemma follows from
\[ \bm{b}(i,j,j,p,q)=[\bm{b}(j,j,i,p), \tau_1^B(A_{p,q})]. \]
$\square$

\vspace{0.5em}

In order to show $\bm{b}(i,j,l,p,j) \in [\mathfrak{b}_n(3),\mathfrak{b}_n(1)]$,
we observe some relations among $\bm{b}(i,j,l,p,q)$s. We use the equality $\equiv$ for the equality in $\mathfrak{b}_n(3)$
modulo $\mathrm{Im}(\tau_4^B)$.
\begin{lem}\label{L-John-4}
For distinct $1 \leq i,j,p,q \leq n$,\\
$(1)$\,\,$\bm{b}(j,i,j,p,q) \equiv \bm{b}(i,j,i,p,q)$.\\
$(2)$\,\,$\bm{b}(i,j,p,q,j) \equiv \bm{b}(i,p,q,i,j)$.\\
$(3)$\,\,$\bm{b}(i,j,p,i,q) \equiv -\bm{b}(j,i,q,j,p)$. \\
$(4)$\,\,$\bm{b}(i,j,p,q,j) \equiv -\bm{b}(i,p,j,q,p)$.
\end{lem}
\textit{Proof.}
By considering the $\mathfrak{S}_n$-orbit, it suffices to show the lemma with specialized $i, j, p, q$.
Recall the element $\bm{b}(1,2,3,4) \in \mathrm{Im}(\tau_3^B)$. Then Part (1), (2) and (3) is obtained from
\[\begin{split}
  \bm{b}(3,2,3,4,1) - \bm{b}(2,3,2,4,1) &= [\bm{b}(1,2,3,4), \tau_1^B(A_{2,3})], \\
  \bm{b}(3,1,2,4,1) - \bm{b}(3,2,4,3,1) &= [\tau_1^B(A_{1,3}), \bm{b}(1,2,3,4)], \\
  \bm{b}(4,3,2,4,1) + \bm{b}(3,4,1,3,2) &= [\bm{b}(1,2,3,4), \tau_1^B(A_{3,4})]. \\
  \end{split}\]
Part (4) is obtained from
\[\begin{split}
   \bm{b}(3,4,1,2,4) &= - \bm{b}(4,3,1,2,4) = \bm{b}(4,3,1,4,2) \overset{(3)}{\equiv} -\bm{b}(3,4,2,3,1) \\
   & \overset{(2)}{\equiv} - \bm{b}(3,1,4,2,1).
  \end{split}\]
This completes the proof of Lemma {\rmfamily \ref{L-John-4}}.
$\square$

\begin{lem}\label{L-John-5}
For distinct $1 \leq i,j,p,q \leq n$, we have
\begin{equation}\begin{split}
 \bm{b}(i,j,& p,q,j) \equiv - \bm{b}(i,p,j,q,p) \equiv -\bm{b}(p,q,i,j,q) \equiv \bm{b}(p,i,q,j,i) \\
  & \equiv \bm{b}(j,q,i,p,q) \equiv - \bm{b}(j,i,q,p,i) \equiv -\bm{b}(q,j,p,i,j) \equiv \bm{b}(q,p,j,i,p) \label{eq-J-1}
\end{split}
\end{equation}
and
\begin{equation}
\begin{split}
 \bm{b}(i,j,& p,i,q) \equiv - \bm{b}(i,q,p,i,j) \equiv -\bm{b}(j,i,q,j,p) \equiv \bm{b}(j,p,q,j,i) \\
  & \equiv -\bm{b}(p,j,i,p,q) \equiv \bm{b}(p,q,i,p,j) \equiv \bm{b}(q,i,j,q,p) \equiv -\bm{b}(q,p,j,q,i). \label{eq-J-2}
\end{split}
\end{equation}
\end{lem}
\textit{Proof.}
We show Equation (\ref{eq-J-1}) since Equation (\ref{eq-J-2}) is obtained by a similar way.
The first equality $\bm{b}(i,j,p,q,j) \equiv - \bm{b}(i,p,j,q,p)$ is Part (4) of Lemma {\rmfamily \ref{L-John-4}}.
We have
\[\begin{split}
  \bm{b}(i,j,p,q,j) & \overset{(2)'}{\equiv} \bm{b}(i,p,q,i,j) \overset{(3)'}{\equiv} -\bm{b}(p,i,j,p,q) 
     \overset{(2)'}{\equiv} -\bm{b}(p,q,i,j,q) \\
   & \overset{(4)'}{\equiv} \bm{b}(p,i,q,j,i).
  \end{split}\]
Here (a)' means Part (a) in Lemma {\rmfamily \ref{L-John-4}}.
By considering the above argument for $\bm{b}(i,p,j,q,p)$ instead of
$\bm{b}(i,j, p,q,j)$, we obtain 
\[ -\bm{b}(i,p,j,q,p) \equiv \bm{b}(j,q,i,p,q) \overset{(4)'}{\equiv} - \bm{b}(j,i,q,p,i). \]
Similarly, by considering the above argument for $\bm{b}(p,i,q,j,i)$ instead of
$\bm{b}(i,j, p,q,j)$, we obtain 
\[ \bm{b}(p,i,q,j,i) \equiv -\bm{b}(q,j,p,i,j) \overset{(4)'}{\equiv} \bm{b}(q,p,j,i,p). \]
Then we obtain all of the equalities.
$\square$

\begin{lem}\label{L-John-6}
For distinct $1 \leq i,j,p,q \leq n$, we have
\[ \bm{b}(i,j,i,p,q) + \bm{b}(j,p,i,q,p)+\bm{b}(i,p,j,q,p) \equiv 0. \]
\end{lem}
\textit{Proof.}
By considering the $\mathfrak{S}_n$-orbit, it suffices to show the lemma with specialized $i, j, p, q$.
Then we see
\[ [\tau_1^B(A_{1,4}), \tau_3^B([A_{1,2},[A_{1,2}, A_{1,3}]])]
   =\bm{b}(2,3,2,1,4) + \bm{b}(3,1,2,4,1)+\bm{b}(2,1,3,4,1).  \]
$\square$

\begin{lem}\label{L-John-9}
For distinct $1 \leq i, j, l, p \leq n$, we have
\[ 2\bm{b}(i,j,l,p,j) \equiv \bm{b}(p,l,l,i,j) + \bm{b}(i,j,j,p,l).  \]
\end{lem}
\textit{Proof.} 
By using Lemma {\rmfamily \ref{L-John-6}} and (\ref{eq-J-1}), we obtain the required equation. $\square$

\vspace{0.5em}

From {\rmfamily \ref{L-John-6}} and {\rmfamily \ref{L-John-9}},
we see $2\bm{b}(i,j,l,p,j) \in [\mathfrak{b}_n(3),\mathfrak{b}_n(1)]$ for $n \geq 4$. In particular,
$\bm{b}(i,j,l,p,j) \in [\mathfrak{b}_n^{\Q}(3),\mathfrak{b}_n^{\Q}(1)]$ for $n \geq 4$.
By combining the above arguments, we have the following.
\begin{pro}\label{P-ab-1}
For any $n \geq 4$, $\mathfrak{b}_n^{\Q}(4)=[\mathfrak{b}_n^{\Q}(3),\mathfrak{b}_n^{\Q}(1)]$, and
$\mathfrak{b}_3(4)=[\mathfrak{b}_3(3),\mathfrak{b}_3(1)]$.
In other words, $H_1(\mathfrak{b}_n, \Q)_4=0$, and $H_1(\mathfrak{b}_n, \Z)_3=0$.
\end{pro}

\subsection{On $H_1(\mathfrak{b}_n,\Z)_5$}\label{Ss-ab-5}

In this subsection we consider $H_1(\mathfrak{b}_n,\Z)_5$. Since $\mathfrak{b}_n(2)=[\mathfrak{b}_n(1),\mathfrak{b}_n(1)]$ as mentioned
in Subsection {\rmfamily \ref{Ss-ab-low}}, we can see
\[ [\mathfrak{b}_n(3),\mathfrak{b}_n(2)] \subset [\mathfrak{b}_n(4),\mathfrak{b}_n(1)] \]
by the Jacobi identity. Hence we have $H_1(\mathfrak{b}_n,\Z)_5=\mathfrak{b}_n(5)/[\mathfrak{b}_n(4),\mathfrak{b}_n(1)]$.
The purpose of this subsection is to show that there exist linearly independent elements in $\mathrm{Ker}(\mathrm{MT}_5^B)$ and
not in $[\mathfrak{b}_n(4),\mathfrak{b}_n(1)]$.

\vspace{0.5em}

For any $1 \leq i, j, l, p, q, r \leq n$, observe
\[\begin{split}
  [x_r,&x_q,x_p,x_l,x_j,x_i] - [x_r,x_q,x_p,x_l,x_i,x_j] - [x_r,x_q,x_p,[x_j,x_i],x_l] \\
   & \hspace{2em} - [x_i,x_j,x_l,[x_q,x_r],x_p] - [x_i,x_j,x_l,x_p,x_r,x_q] + [x_i,x_j,x_l,x_p,x_q,x_r] \\
   &= [x_r,x_q,x_p,x_l,[x_j,x_i]] - [x_r,x_q,x_p,[x_j,x_i],x_l] \\
   & \hspace{2em} -[x_i,x_j,x_l,[x_q,x_r],x_p] - [x_i,x_j,x_l,x_p,[x_r,x_q]] \\
   &= [x_r,x_q,x_p,[x_l,[x_j,x_i]]] + [x_i,x_j,x_l,[[x_r,x_q],x_p]] = 0,
  \end{split}\]
and set
\[\begin{split}
   \bm{b}(i,j,l,p,q,r) := & x_i^* \otimes [x_r,x_q,x_p,x_l,x_j,x_i] - x_j^* \otimes [x_r,x_q,x_p,x_l,x_i,x_j] \\
   & - x_l^* \otimes [x_r,x_q,x_p,[x_j,x_i],x_l] - x_p^* \otimes [x_i,x_j,x_l,[x_q,x_r],x_p] \\
   & - x_q^* \otimes [x_i,x_j,x_l,x_p,x_r,x_q] + x_r^* \otimes [x_i,x_j,x_l,x_p,x_q,x_r] \in \mathfrak{b}_n(5).
  \end{split}\]
We remark that $\bm{b}(i,j,l,p,q,r)$ corresponds to an internal trivalent tree
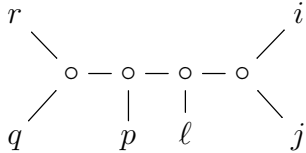
\[ \xymatrix@R=3mm@C=3mm{
r\ar@{-}[dr]&&&&&i \cr
& \circ\ar@{-}[r] & \circ\ar@{-}[r] &\circ\ar@{-}[r] &\circ\ar@{-}[ru]\ar@{-}[rd] &\cr
q\ar@{-}[ru]&&p\ar@{-}[u]&\ell\ar@{-}[u]&&j
} \]
\begin{lem}\label{L-ab-5-1}
For any $1 \leq i, j, l, p, q, r \leq n$,
\[\begin{split}
  \mathrm{MT}_5^B( & \bm{b}(i,j,l,p,q,r)) \\
   &= 2(\delta_{ir} \bm{e}_{x_i,x_j,x_l,x_p,x_q} - \delta_{iq} \bm{e}_{x_j,x_l,x_p,x_q,x_r} - \delta_{j,r} \bm{e}_{x_i,x_j,x_l,x_p,x_q}
        + \delta_{j,q} \bm{e}_{x_i,x_l,x_p,x_q,x_r}) 
  \end{split}\]
\end{lem}
\textit{Proof.}
This lemma follows from Lemma {\rmfamily \ref{L-S03-1}}, by a straightforward calculation. $\square$

\vspace{0.5em}

For any $1 \leq i, j, l, p \leq n$, consider
\[\begin{split}
   \bm{n}(i,j,l) &:=\bm{b}(l,j,i,l,j,i)-\bm{b}(l,j,l,i,j,i) \in \mathfrak{b}_n(5,(2^3)), \\
   \bm{n}(i,j,l,p) &:=\bm{b}(p,j,j,i,l,i) \in \mathfrak{b}_n(5,(2^2,1^2)).
  \end{split}\]
From Lemma {\rmfamily \ref{L-ab-5-1}}, we see that $\bm{n}(i,j,l), \bm{n}(i,j,l,p) \in \mathrm{Ker}(\mathrm{MT}_5^B)$.
We show that these elements are non-trivial and do not belong to $[\mathfrak{b}_n(4),\mathfrak{b}_n(1)]$.
First we show that these elements are non-trivial.

\begin{lem}\label{L-ab-5-2}
For any $1 \leq i, j, l, p \leq n$, $\bm{n}(i,j,l), \bm{n}(i,j,l,p) \neq 0 \in \mathfrak{b}_n(5)$.
\end{lem}
\textit{Proof.}
First we consider $\bm{n}(i,j,l)$s.
By observing the $\mathfrak{S}_n$-orbit, it suffices to show the lemma with specialized $i, j, l$. So we show $\bm{n}(1,2,3) \neq 0 \in \mathfrak{b}_n(5)$.
Then we have
\[\begin{split}
  \bm{b}(3,2,1,3,2,1) & = x_3^* \otimes [x_1,x_2,x_3,x_1,x_2,x_3] - x_2^* \otimes [x_1,x_2,x_3,x_1,x_3,x_2] \\
                      & \hspace{2em} +x_1^* \otimes [x_1,x_2,x_3,[x_3,x_2],x_1] - x_3^* \otimes [x_3,x_2,x_1,[x_2,x_1],x_3] \\
                      & \hspace{2em} -x_2^* \otimes [x_3,x_2,x_1,x_3,x_1,x_2] +x_1^* \otimes [x_3,x_2,x_1,x_3,x_2,x_1], \\
  \bm{b}(3,2,3,1,2,1) &= x_3^* \otimes [x_1,x_2,x_1,x_3,x_2,x_3] - x_2^* \otimes [x_1,x_2,x_1,x_3,x_3,x_2] \\
                      & \hspace{2em} +x_3^* \otimes [x_1,x_2,x_1,[x_3,x_2],x_3] - x_1^* \otimes [x_3,x_2,x_3,[x_2,x_1],x_1] \\
                      & \hspace{2em} -x_2^* \otimes [x_3,x_2,x_3,x_1,x_1,x_2] +x_1^* \otimes [x_3,x_2,x_3,x_1,x_2,x_1].
  \end{split}\]
By considering the isomorphism $H \otimes_{\Z} \mathcal{L}_n(k) \cong \mathfrak{p}_n(k)$ given $x_i \otimes X \mapsto x_i^* \otimes [X, x_i]$,
we observe the part $x_3^* \otimes [-, x_3]$ in $\bm{n}(1,2,3)$.
For example, by rewriting $[x_1,x_2,x_3,x_1,x_2]$ as a linear combination of the basic commutators of weight five, we have
\[\begin{split}
   [x_1,x_2,x_3,x_1,x_2] =[ & x_3,x_1,x_2,[x_2,x_1]] - [x_2,x_1,x_2,[x_3,x_1]] - [x_2,x_1,x_1,[x_3,x_2]] \\
    & - [x_2,x_1,x_1,x_2,x_3]. 
  \end{split}\]
Similarly, by rewriting $[x_3,x_2,x_1,[x_2,x_1]]$, $[x_1,x_2,x_1,x_3,x_2]$ and $[x_1,x_2,x_1,[x_3,x_2]]$ as
a linear combination of the basic commutators of weight five. we can see that
$x_3^* \otimes [x_2,x_1,x_2,[x_3,x_1], x_3]$ appears only from $x_3^* \otimes [x_1,x_2,x_3,x_1,x_2,x_3]$, and hence $\bm{n}(1,2,3) \neq 0 \in \mathfrak{b}_n(5)$.

Next we consider $\bm{n}(i,j,l,p)$. By a similar argument as above, it suffices to show that $\bm{n}(1,2,3,4) \neq 0 \in \mathfrak{b}_n(5)$.
By observing the term $x_4^* \otimes [-,x_4]$ in $\bm{n}(1,2,3,4)$, we can immediately see it. $\square$

\vspace{0.5em}

Now we show $\bm{n}(i,j,l) \notin [\mathfrak{b}_n(4), \mathfrak{b}_n(1)]$. By considering the $\mathfrak{S}_n$-orbit, it suffices to show
$\bm{n}(1,2,3) \notin [\mathfrak{b}_n(4), \mathfrak{b}_n(1)]$.
In general, the standard inclusion $F_3 \rightarrow F_n$ defined by $x_i \mapsto x_i$ for $1 \leq i \leq 3$ induces an injective Lie algebra
homomorphism $\mathcal{L}_3 \rightarrow \mathcal{L}_n$. This injection also induces an injective Lie algebra
homomorphism $\mathfrak{p}_3 \rightarrow \mathfrak{p}_n$. Through this injection, we consider $\mathfrak{p}_3$ as a Lie subalgebra of $\mathfrak{p}_n$.
Similarly, we consider $\mathfrak{b}_3$ as a Lie subalgebra of $\mathfrak{b}_n$.
From this viewpoint, we have $\mathfrak{b}_3(k) \subset \mathfrak{b}_n(k)$ for each $k \geq 1$.
\begin{lem}\label{L-ab-5-2}
If $\bm{n}(1,2,3) \in [\mathfrak{b}_n(4), \mathfrak{b}_n(1)]$, then $\bm{n}(1,2,3) \in [\mathfrak{b}_3(4), \mathfrak{b}_3(1)]$.
\end{lem}
\textit{Proof.}
In general, for any $k \geq 1$
the Lie bracket of $x_{i_{k+1}}^* \otimes [x_{i_1}, \ldots, x_{i_{k+1}}] \in \mathfrak{p}_n(k)$ and
$x_{j_2}^* \otimes [x_{j_1},x_{j_2}] \in \mathfrak{p}_n(1)$ in $\mathfrak{p}_n(k+1)$ is calculated as
\[\begin{split}
  [x_{i_{k+1}}^* & \otimes [x_{i_1}, \ldots, x_{i_{k+1}}], x_{j_2}^* \otimes [x_{j_1},x_{j_2}]] \\
   &= \sum_{m=1}^{k+1} \delta_{i_m, j_2} x_{i_{k+1}}^* \otimes [x_{i_1}, \ldots, [x_{j_1},x_{j_2}], \ldots, x_{i_{k+1}}] \\
   & \hspace{1em} - \delta_{j_1, i_{k+1}} x_{j_2}^* \otimes [x_{i_1}, \ldots, x_{i_{k+1}},x_{j_2}]
     - \delta_{j_2, i_{k+1}} x_{j_2}^* \otimes [x_{j_1},[x_{i_1}, \ldots, x_{i_{k+1}}]].
  \end{split}\]
Thus, any component in the commutators in $\mathcal{L}_n(k+2)$ appearing the right hand side of the above equation is in
$\{ x_{i_1}, \ldots, x_{i_{k+1}}, x_{j_1},x_{j_2} \}$. Moreover, any element in $\{ x_{i_1}, \ldots, x_{i_{k+1}}, x_{j_1},x_{j_2} \}$
appears at least one in the components of the commutators in appearing the right hand side of the above equation.

On the other hand, if $\{ a_1, \ldots, a_{k+1} \}$ and $\{ b_1, \ldots, b_{k+1} \}$ for $1 \leq a_m, b_m \leq n$
are different set with multiplicity, and $[a_1, \ldots, a_{k+1}], [b_1, \ldots, b_{k+1}] \neq 0 \in \mathcal{L}_n(k+1)$,
then we can see that $[a_1, \ldots, a_{k+1}]$ and $[b_1, \ldots, b_{k+1}]$ are linearly independent in $\mathcal{L}_n(k+1)$
by rewriting them as a linear combination of Hall's basic commutators. Hence, we see
$x_{a_{k+1}}^* \otimes [a_1, \ldots, a_{k+1}]$ and $x_{b_{k+1}}^* \otimes [b_1, \ldots, b_{k+1}]$ are also linearly independent in $\mathfrak{p}_n$.

Therefore, by observing $\mathfrak{b}_n \subset \mathfrak{p}_n$, we see that
if $\bm{n}(1,2,3) \in [\mathfrak{b}_n(4), \mathfrak{b}_n(1)]$ then $\bm{n}(1,2,3) \in [\mathfrak{b}_3(4), \mathfrak{b}_3(1)]$.
$\square$

\begin{lem}\label{L-ab-5-3}
$\bm{n}(1,2,3) \notin [\mathfrak{b}_3(4), \mathfrak{b}_3(1)]$.
\end{lem}
\textit{Proof.}
Assume $\bm{n}(1,2,3) \in [\mathfrak{b}_3(4), \mathfrak{b}_3(1)]$.
Since $\bm{n}(1,2,3) \in \mathfrak{S}_n \cdot \mathfrak{b}_3(5, (2^3))$, there is the only one possibility that
$\bm{n}(1,2,3) \in [\mathfrak{S}_3 \cdot \mathfrak{b}_3(4, (2^2,1)), \mathfrak{b}_3(1)]$.
From a result in Subsection {\rmfamily \ref{Sss-4}}, we have
\[ \mathfrak{S}_3 \cdot \mathfrak{b}_3(4, (2^2,1)) = \langle \bm{b}(2,1,1,2,3), \bm{b}(3,1,1,3,2), \bm{b}(3,2,2,3,1) \rangle_{\Z} \]
and $\mathfrak{b}_3(1)$ is generated by $\tau_1^B(A_{1,2})$, $\tau_1^B(A_{1,3})$ and $\tau_1^B(A_{2,3})$.
The Lie brackets of the above generators of $\mathfrak{S}_3 \cdot \mathfrak{b}_3(4, (2^2,1))$ and $\mathfrak{b}_3(1)$
that have terms in $\mathfrak{S}_n \cdot \mathfrak{b}_3(5, (2^3))$ are
\[\begin{split}
  Z_1 &:= [\bm{b}(2,1,1,2,3), \tau_1^B(A_{1,3})], \hspace{1em} Z_2 := [\bm{b}(2,1,1,2,3), \tau_1^B(A_{2,3})], \\
  Z_3 &:= [\bm{b}(3,1,1,3,2), \tau_1^B(A_{1,2})], \hspace{1em} Z_4 := [\bm{b}(3,1,1,3,2), \tau_1^B(A_{2,3})], \\
  Z_5 &:= [\bm{b}(3,2,2,3,1), \tau_1^B(A_{1,2})], \hspace{1em} Z_6 := [\bm{b}(3,2,2,3,1), \tau_1^B(A_{1,3})].
  \end{split}\]
By using the formula (\ref{eq-der-cal}), we can calculate the above elements, and see that each of $Z_i$ has terms in 
$\mathfrak{S}_n \cdot \mathfrak{b}_3(5, (2^3))$ and $\mathfrak{S}_n \cdot \mathfrak{b}_3(5, (3,2,1))$.
Write the $\mathfrak{S}_n \cdot \mathfrak{b}_3(5, (3,2,1))$ part of $Z_i$ as $\overline{Z_i}$.
Since $\bm{n}(1,2,3) \in [\mathfrak{S}_3 \cdot \mathfrak{b}_3(4, (2^2,1)), \mathfrak{b}_3(1)]$,
$\bm{n}(1,2,3)$ is written as $a_1Z_1 + \cdots + a_6Z_6$. Then we see $a_1\overline{Z_1} + \cdots + a_6\overline{Z_6}=0$.

\vspace{0.5em}

On the other hand, we have
{\small
\[\begin{split}
  \overline{Z_1} = & -x_1^* \otimes [x_3,x_1,x_2,x_1,x_2,x_1] + x_1^* \otimes [x_3,x_1,x_2,[x_2,x_1],x_1] -x_1^* \otimes [x_2,x_1,x_1,x_2,x_3,x_1] \\
   & + x_2^* \otimes [x_2,x_1,x_1,[x_3,x_1],x_2] + x_2^* \otimes [x_3,x_1,x_2,x_1,x_1,x_2] +x_3^* \otimes [x_2,x_1,x_1,x_2,x_1,x_3], \\
  \overline{Z_2} = & -x_1^* \otimes [x_1,x_2,x_2,[x_3,x_2],x_1] - x_1^* \otimes [x_3,x_2,x_1,x_2,x_2,x_1] +x_2^* \otimes [x_3,x_2,x_1,x_2,x_1,x_2] \\
   & - x_2^* \otimes [x_3,x_2,x_1,[x_1,x_2],x_2] + x_2^* \otimes [x_1,x_2,x_2,x_1,x_3,x_2] -x_3^* \otimes [x_1,x_2,x_2,x_1,x_2,x_3], \\
  \overline{Z_3} = & -x_1^* \otimes [x_2,x_1,x_3,x_1,x_3,x_1] + x_1^* \otimes [x_2,x_1,x_3,[x_3,x_1],x_1] -x_1^* \otimes [x_3,x_1,x_1,x_3,x_2,x_1] \\
   & + x_2^* \otimes [x_3,x_1,x_1,x_3,x_1,x_2] + x_3^* \otimes [x_3,x_1,x_1,[x_2,x_1],x_3] +x_3^* \otimes [x_2,x_1,x_3,x_1,x_1,x_3], \\
  \overline{Z_4} = & -x_1^* \otimes [x_1,x_3,x_3,[x_2,x_3],x_1] + x_1^* \otimes [x_2,x_3,x_1,x_3,x_3,x_1] -x_2^* \otimes [x_1,x_3,x_3,x_1,x_3,x_2] \\
   & + x_3^* \otimes [x_2,x_3,x_1,x_3,x_1,x_3] - x_3^* \otimes [x_2,x_3,x_1,[x_1,x_3],x_3] +x_3^* \otimes [x_1,x_3,x_3,x_1,x_2,x_3], \\
  \overline{Z_5} = & \hspace{1em} x_1^* \otimes [x_3,x_2,x_2,x_3,x_2,x_1] - x_2^* \otimes [x_1,x_2,x_3,x_2,x_3,x_2] +x_2^* \otimes [x_1,x_2,x_3,[x_3,x_2],x_2] \\
   & - x_2^* \otimes [x_3,x_2,x_2,x_3,x_1,x_2] + x_3^* \otimes [x_3,x_2,x_2,[x_1,x_2],x_3] +x_3^* \otimes [x_1,x_2,x_3,x_2,x_2,x_3], \\
  \overline{Z_6} = & -x_1^* \otimes [x_2,x_3,x_3,x_2,x_3,x_1] - x_2^* \otimes [x_2,x_3,x_3,[x_1,x_3],x_2] -x_2^* \otimes [x_1,x_3,x_2,x_3,x_3,x_2] \\
   & + x_3^* \otimes [x_1,x_3,x_2,x_3,x_2,x_3] - x_3^* \otimes [x_1,x_3,x_2,[x_2,x_3],x_3] +x_3^* \otimes [x_2,x_3,x_3,x_2,x_1,x_3]. \\
 \end{split}\]
}
\noindent
We remark that since $Z_2=- [\bm{b}(1,2,2,1,3), \tau_1^B(A_{2,3})] = - \sigma \cdot Z_1$ for $\sigma=(1,2) \in \mathfrak{S}_3$, we see
$\overline{Z_2}=-\sigma \cdot \overline{Z_1}$, Similarly we can calculate $\overline{Z_i}$ for any $2 \leq i \leq 6$ from $\overline{Z_1}$.
Observe the $x_3^* \otimes [-, x_3]$ part in $a_1\overline{Z_1} + \cdots + a_6\overline{Z_6}=0$.
An element type of $x_3^* \otimes [C, x_3]$ where $C$ is a commutator of weight $5$, and whose component consists of three $x_1$s and two $x_2$s only appears in
$\overline{Z_1}$ as $x_3^* \otimes [x_2,x_1,x_1,x_2,x_1,x_3]$. This means $a_1=0$.
By the same argument, by observing elements type of $x_3^* \otimes [C, x_3]$ where $C$ is a commutator of weight $5$, and whose component consists of two $x_1$s and three $x_2$s,
we obtain $a_2=0$. Similarly, by observing the $x_2^* \otimes [-, x_2]$ part, we obtain $a_3=0$ and $a_4=0$.
Finally, by observing the $x_1^* \otimes [-, x_1]$ part, we obtain $a_5=0$ and $a_6=0$.
This is a contradiction since $\bm{n}(1,2,3)$ is non-trivial by Lemma {\rmfamily \ref{L-ab-5-2}}.
$\square$

\vspace{0.5em}

From Lemmas {\rmfamily \ref{L-ab-5-2} and {\rmfamily \ref{L-ab-5-3}}, we obtain the following proposition.

\begin{pro}\label{P-ab-5-4}
For any $n \geq 3$, we have
$\bm{n}(1,2,3) \notin [\mathfrak{b}_n(4), \mathfrak{b}_n(1)]$.
\end{pro}

As a corollary to Proposition {\rmfamily \ref{P-ab-5-4}}, we obtain the following consequece.
\begin{thm}\label{T-ninj-MT}
For $n \geq 3$, the surjective Lie algebra homomorphism
\[ \overline{\Theta|_{\mathfrak{b}_n}} : H_1(\mathfrak{b}_n, \Z) \rightarrow \mathrm{Im}(\tau_1) \oplus
   \bigoplus_{m \geq 1} \mathrm{Im}(\mathrm{MT}_{2m+1}^B)  \]
induced from $\Theta|_{\mathfrak{b}_n}$ is not injective.
\end{thm}

\begin{rem}\label{R-ab-5-5}
By using a computer calculation, we can show that 
\[
\bm{n}(1,2,3,4) \notin [\mathfrak{b}_n(4), \mathfrak{b}_n(1)].
\]
The method is the same as above. Namely, we can check that  
\[
\bm{n}(1,2,3,4) \notin [\mathfrak{S}_4 \cdot \mathfrak{b}_4(4, (2^2,1)), 
   \mathfrak{b}_4(1)] \oplus [\mathfrak{S}_4 \cdot \mathfrak{b}_4(4, (2,1^3)), \mathfrak{b}_4(1)]. 
   \]
\quad Based on these results, we expect that 
\[ \mathrm{dim}_{\Q}(\mathrm{Ker}(\mathrm{MT}_{5,\Q}^B)/[\mathfrak{b}_n^{\Q}(4), \mathfrak{b}_n^{\Q}(1)]) =\binom{n}{3}+\binom{n}{4} \]
with a basis
\[ \{ \bm{n}(i_1,i_2,i_3) \,|\, 1 \leq i_1<i_2<i_3 \leq n \} \cup \{ \bm{n}(i_1,i_2,i_3,i_4)  \,|\, 1 \leq i_1<i_2<i_3<i_4 \leq n \}. \]
That is true for $n \leq 7$. 
\end{rem}

\section{Acknowledgments}

\vspace{0.5em}

We would like to thank Takuya Sakasai, Yusuke Kuno, Masatoshi Sato for various comments for the
the special derivation Lie algebras of free Lie algebras.
In particular, we would like to thank Yusuke Kuno for communicating the non-finite generation of $H_1(\mathfrak{b}_2,\Z)$ with Soul\'{e} elements.
We also would like to thank the referee for his/her valuable comments for the original version of the paper.

The second author is supported by JSPS KAKENHI Grant Number 25K06988.

\end{document}